\documentclass[times,review,10pt]{elsarticle}

\usepackage{amssymb}
\usepackage{amsmath,amssymb,amsfonts}
\usepackage{cases}
\usepackage{color,xcolor}
\usepackage{graphicx}
\usepackage{subfigure}
\usepackage{algorithm}
\usepackage{algorithmic}
\usepackage{epstopdf}
\usepackage{multirow}
\usepackage{rotating}
\usepackage{booktabs}
\usepackage{tablefootnote}
\usepackage{tikz}
\usepackage{bm}
\usepackage{color,xcolor}
\usepackage[colorlinks, linkcolor=blue, anchorcolor=blue, citecolor=blue]{hyperref}

\newtheorem{Theorem}{Theorem}[section]

\newtheorem{Remark}[Theorem]{Remark}

\journal{Elsevier}

\begin{document}

\begin{frontmatter}

\title{Enhancing Unsupervised Feature Selection via Double Sparsity Constrained Optimization}

\author[1]{Xianchao Xiu}
\ead{xcxiu@shu.edu.cn}

\author[1]{Anning Yang}
\ead{yanganning1029@163.com}

\author[1]{Chenyi Huang}
\ead{huangchenyi@shu.edu.cn}

\author[2]{Xinrong Li \corref{cor1}}
\ead{lixinrong@mail.neu.edu.cn}

\author[3]{Wanquan Liu}
\ead{liuwq63@mail.sysu.edu.cn}

\address[1]{School of Mechatronic Engineering and Automation, Shanghai University, Shanghai 200444, China}
\address[2]{National Frontiers Science Center for Industrial Intelligence and Systems Optimization, Northeastern University, Shenyang 110819, China}
\address[3]{School of Intelligent Systems Engineering, Sun Yat-sen University, Guangzhou 528406, China}

\cortext[cor1]{Corresponding author}

\begin{abstract}
Unsupervised feature selection (UFS) is widely applied in machine learning and pattern recognition. However, most of the existing methods only consider a single sparsity, which makes it difficult to select valuable and discriminative feature subsets from the original high-dimensional feature set. In this paper, we propose a new UFS method called DSCOFS via embedding double sparsity constrained optimization into the classical principal component analysis (PCA) framework. Double sparsity refers to using $\ell_{2,0}$-norm and  $\ell_0$-norm to simultaneously constrain variables, by adding the sparsity of different types, to achieve the purpose of improving the accuracy of identifying differential features. The core is that $\ell_{2,0}$-norm can remove  irrelevant and redundant features, while $\ell_0$-norm can filter out irregular noisy features, thereby complementing $\ell_{2,0}$-norm to improve discrimination. An effective proximal alternating minimization method is proposed to solve the resulting nonconvex nonsmooth model. Theoretically, we rigorously prove that the sequence generated by our method globally converges to a stationary point. Numerical experiments on three synthetic datasets and eight real-world datasets demonstrate the effectiveness, stability, and convergence of the proposed method. In particular, the average clustering accuracy (ACC) and normalized mutual information (NMI) are improved by at least 3.34\% and 3.02\%, respectively, compared with the state-of-the-art methods. More importantly, two common statistical tests and a new feature similarity metric verify the advantages of double sparsity. All results suggest that our proposed DSCOFS provides a new perspective for  feature selection.
\end{abstract}


\end{frontmatter}

\section{Introduction}
As a popular dimensionality reduction method, feature selection has become an important technology for machine learning,  pattern recognition, and image processing. Feature selection aims to eliminate redundant, irrelevant, and noisy features, thereby reducing training time and enhancing interpretability \cite{fan2024learning}. 
However, real-world high-dimensional data often lack labels \cite{zhou2024unsupervised}. Fortunately, unsupervised feature selection (UFS) has emerged and attracted considerable attention due to its capacity to select the most representative features without relying on labels. Nowadays, UFS has been successfully applied in many fields ranging from image classification \cite{huang2023robust} to remote sensing \cite{uddin2021information}, wireless communications \cite{liu2022deep}, and genetic engineering \cite{saberi2022dual}. The interested readers can refer to the excellent surveys \cite{dhal2022comprehensive,li2024exploring}.

During the past few decades, a large number of UFS methods based on spectral analysis have been proposed and studied, among which the representatives are but not limited to Laplacian score (LapScore) \cite{he2005laplacian}, unsupervised discriminative feature selection (UDFS) \cite{yang2011l}, structured optimal graph feature selection (SOGFS) \cite{nie2019structured}, and robust neighborhood embedding (RNE) \cite{liu2020robust}. However, the performance of these spectral-based methods tends to be affected by the fixed graphs constructed in high-dimensional spaces, which makes them sensitive to redundant and noisy features, and the computational cost increases greatly as the number of samples increases.

In order to address the above problems, principal component analysis (PCA) and sparse PCA (SPCA) \cite{zou2006sparse} have been increasingly studied in the field of UFS. For example, Chang \textit{et al.} \cite{chang2016convex} constructed an effective convex SPCA model by introducing $\ell_{2,1}$-norm to measure the loss and regularization. Yi \textit{et al.} \cite{yi2019adaptive} proposed adaptive weighted SPCA to select significant features from noisy data robustly. Recently, Li \textit{et al.} \cite{li2021sparse} considered  $\ell_{2,p}$-norm $(0<p\leq 1)$ regularized SPCA and proposed SPCAFS with an alternating minimization optimization algorithm. In addition, Zheng \textit{et al.} \cite{zheng2023fast} used the reconstruction matrix to reformulate SPCA as a positive semidefinite (PSD) problem. It is worth noting that these SPCA-based UFS methods only consider convex or nonconvex relaxations of $\ell_{2,0}$-norm. As discussed in \cite{pang2018efficient,wang2023sparse}, this may lead to insufficient sparse representation. Motivated by this, Chen \textit{et al.} \cite{chen2022fast} embedded an adaptive bipartite graph in the subspace using the original data and the selected anchors, and then applied an $\ell_{2,0}$-norm constraint to the projection matrix for feature selection. Nie \textit{et al.} \cite{nie2022learning} studied the feature-sparsity constrained PCA (FSPCA) problem with $\ell_{2,0}$-norm constraints based on the rank of the data covariance matrix.

The $\ell_{2,0}$-norm and its relaxed forms only represent a type of structural sparsity \cite{cui2021fused}. They can only handle the irrelevant and redundancy between features. Nevertheless, when the data exhibits element-wise sparsity, as indicated by the $\ell_{0}$-norm,  it may lead to poor feature selection and, in turn, bad clustering results \cite{liu2024towards}. Although some studies have applied double sparsity to signal recovery \cite{bian2023joint}, radar imaging \cite{zhang2021micro}, compressive sensing \cite{zhou2022computing}, and feature selection \cite{guo2023double}, they are not constrained to the same variable. It is worth noting that there have been instances of successful research in the brain imaging predictor identification \cite{wang2011sparse} and gene expression \cite{hu2019dstpca}, where the incorporation of $\ell_{2,1}$-norm and $\ell_1$-norm regularization into the objective function for the same variable has been employed. In fact, the $\ell_{2,0}$-norm and $\ell_0$-norm constraints are more flexible and accurate than regularization and relaxation methods in determining sparsity \cite{sun2023gradient}, allowing direct access to the desired features without the need to compute the transformation matrix score \cite{nie2024row}. \textit{Therefore, a natural question is whether it is possible to integrate the original double sparsity constraints, i.e., $\ell_{2,0}$-norm plus $\ell_{0}$-norm, into a PCA framework for better feature selection?}

\begin{figure}[t]
\centering
    \hspace{-4.5mm}
	\subfigcapskip=-5pt
	\subfigure[Original]{
		\centering
        \label{a}
        \includegraphics[width=0.32\textwidth]{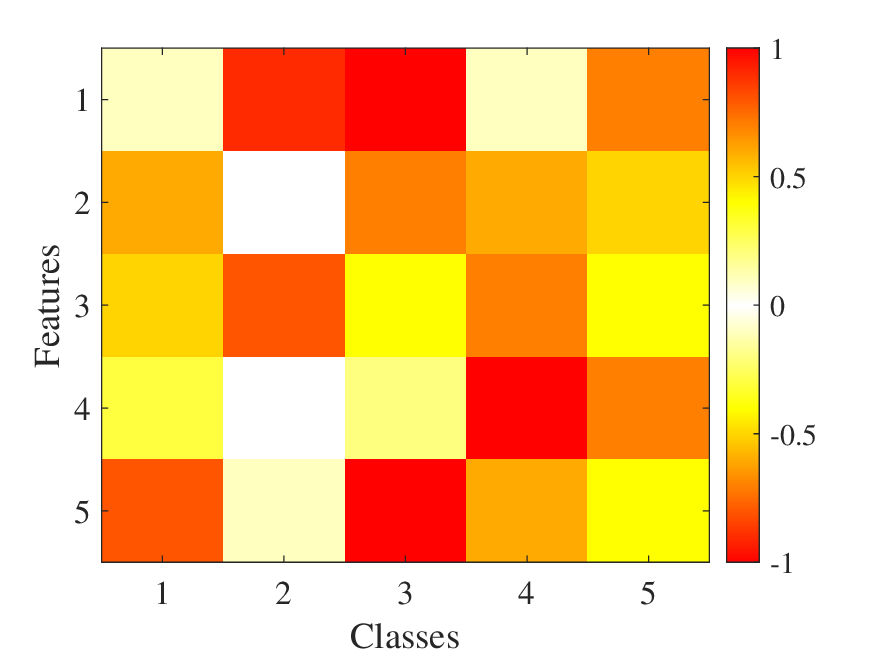}
	}
    \hspace{-2mm} 
    \subfigcapskip=-5pt
	\subfigure[Obtained by $\ell_{2,0}$-norm]{		
		\label{b}
		\centering
        \includegraphics[width=0.32\textwidth]{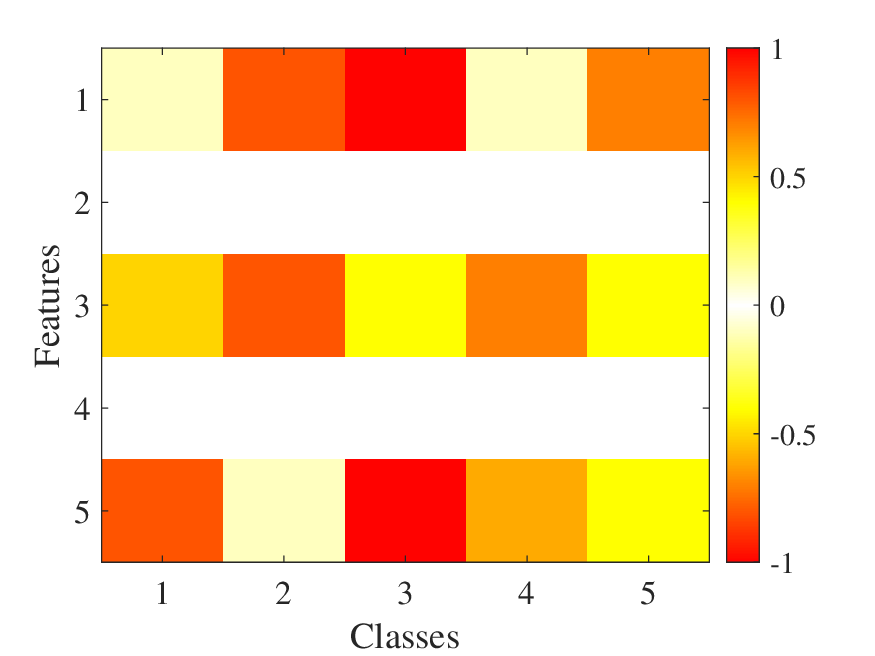}				
	}
    \hspace{-4.5mm}

    \vspace{-2.5mm}  
    \hspace{-4.5mm}
	\subfigcapskip=-5pt
	\subfigure[Obtained by $\ell_0$-norm]{
		\centering
        \label{a}
        \includegraphics[width=0.32\textwidth]{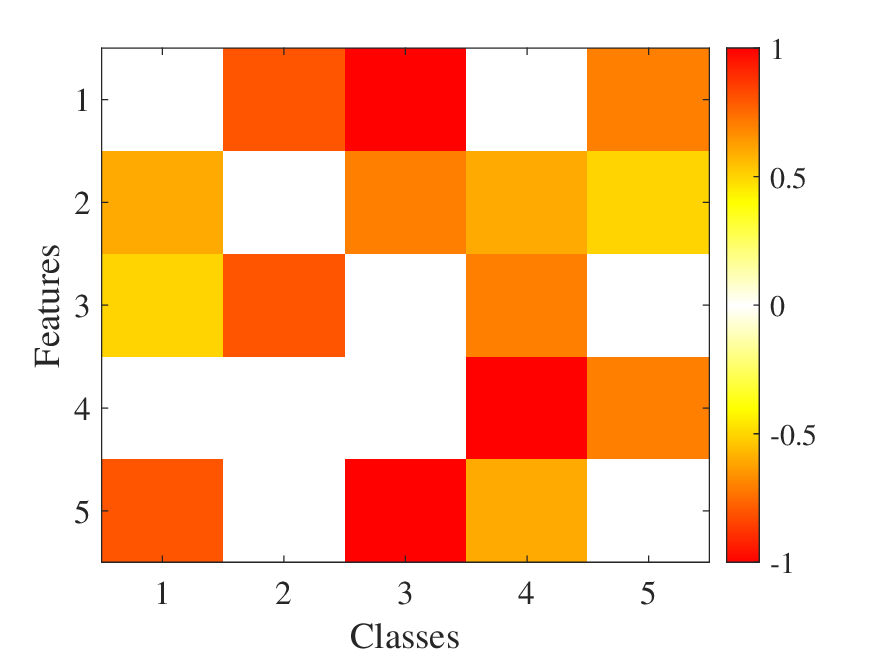}
	}
    \hspace{-2mm} 
    \subfigcapskip=-5pt
	\subfigure[Obtained by  double sparsity]{		
		\label{b}
		\centering
        \includegraphics[width=0.32\textwidth]{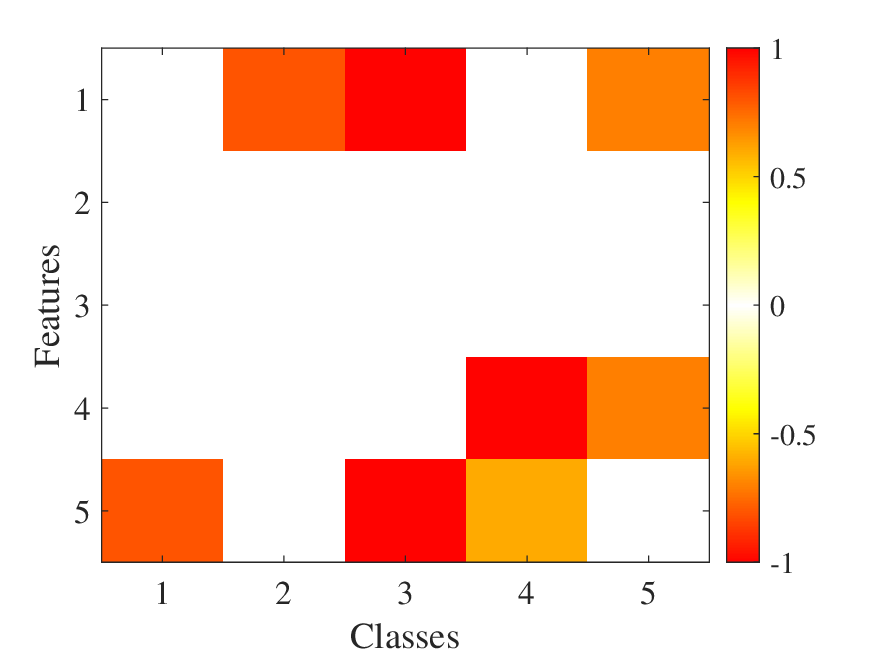}				
	}
    \hspace{-4.5mm}
    \vspace{-0.2cm}
\caption{The examples of results obtained by different sparsity constraints, where white means the value of the element is 0.}
\centering
\label{plot-heat}
\end{figure}

In this paper, we propose a novel UFS method called DSCOFS with the help of double sparsity constrained optimization, where $\ell_{2,0}$-norm ensures the global structural sparsity and has better interpretability, while $\ell_0$-norm takes the local individual sparsity of data elements into account. As shown in Figure \ref{plot-heat}, double sparsity can not only remove irrelevant and redundant features but also select more representative features. In fact, double sparsity can be viewed as an attempt to bridge the gap between single sparsity and actual sparsity. Experimental results show that compared with the state-of-the-art UFS methods, the average clustering accuracy (ACC) and normalized mutual information (NMI) of all datasets are improved at least by at least 3.34\% and 3.02\%, respectively. Figure \ref{plot-time} indicates although DSCOFS takes a little more time than traditional methods, there is a significant improvement in terms of ACC.

\begin{figure}[t] \label{plot-time}
\centering
\includegraphics[width=0.70\textwidth]{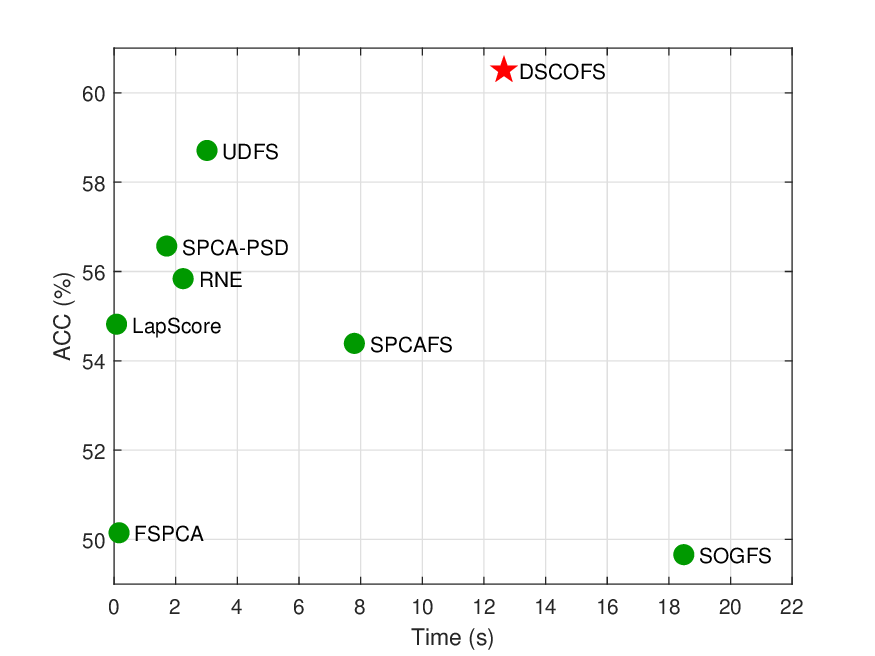}
\vspace{-0.2cm}
\caption{The ACC performance of our proposed DSCOFS and other UFS methods, including LapScore, UDFS, SOGFS, RNE, FSPCA, SPCAFS, SPCA-PSD, under time on the Isolet dataset.}
\label{plot-time}
\end{figure}

The main contributions of this paper can be summarized in four aspects.

\begin{itemize}
\item We establish a novel double sparsity constrained PCA with $\ell_{2,0}$-norm and $\ell_0$-norm for UFS, which has not been sufficiently investigated before. Unlike the existing single sparsity constrained and regularized PCA, double sparsity is able to obtain more discriminative features and make our proposed DSCOFS perform better when facing complicated data with different structural features.
\item We develop a proximal alternating minimization algorithm that combines an exact penalty function method with the hard threshold operators to solve the model. All the resulting subproblems admit closed-form solutions or can be computed by fast solvers, thus providing an efficient optimization strategy.
\item We prove that after a finite number of iterations, the objective function of the developed algorithm is strictly nonincreasing, which is similar to the conclusion in most literature. We further rigorously show that the generated sequence globally converges to a stationary point of our model without any further assumption.
\item We define a new metric called feature similarity ratio (FSR). As shown in Table \ref{Ablation}, the features extracted by our proposed DSCOFS  are different from those extracted by the single $\ell_{2,0}$-norm, and even only 52\% similarity on the Isolet dataset. In addition, the Friedman test and post-hoc test validate that our proposed DSCOFS is significantly different from other comparison methods.
\end{itemize}


The rest of this paper is organized as follows. Section \ref{pre} introduces notations and related works. Section \ref{method} gives our mathematical model with the optimization algorithm and convergence analysis. Section \ref{expriments} evaluates the effectiveness and superiority of our proposed method. Section \ref{conclusion} presents the conclusion and future work. The detailed proof is provided in the appendix.

\section{Preliminaries}\label{pre}

This section first provides some necessary notations and then describes related works about SPCA.

\subsection{Notations}

Throughout this paper, matrices are represented by uppercase letters such as $X$.  For any matrix $X\in\mathbb{R}^{d\times m}$, the $ij$-th, $i$-th row, $j$-th column are denoted by $X_{ij}$, $X_{i:}$, $X_{:j}$, respectively. The Euclidean inner product of two matrices $X, Y\in\mathbb{R}^{d\times m}$ is defined as $\langle X, Y\rangle=\textrm{Tr}(X^\top Y)$, where $\textrm{Tr}(\cdot)$ indicates the trace and $X^\top$ represents the transpose of $X$. The $\ell_{0}$-norm and $\ell_{2,0}$-norm are, respectively, defined as $\|X\|_{0}=\textrm{card}(\{ij \mid |X_{ij}|\neq 0\})$ and $\|X\|_{2,0}=\textrm{card}(\{i\mid \|X_{i:}\|_2\neq 0\})$, where $\textrm{card}(\cdot)$ denotes the cardinality. Let $\ell_{2,p}$-norm ($0<p\leq 1$) be $\|X\|_{2,p}^p=\sum_{i=1}^d \|X_{i:}\|_2^p$. When $p=1$, it is called $\ell_{2,1}$-norm. Besides, $\|X\|_\textrm{F}$ represents the Frobenius norm, and $I_m$ is the $m\times m$ identity matrix. Further notation will be introduced wherever it occurs.

\subsection{SPCA Basics}

Given a data matrix $A = [A_{:1}, A_{:2}, A_{:3}, \cdots, A_{:n}] \in \mathbb{R}^{d\times n}$, where each column represents a sample with $d$ features. Without loss of generality, assume that the data has been centered by $\sum_{i=1}^n A_{:j} =0$. The purpose of PCA is to find a low-dimensional transformation matrix $X\in \mathbb{R}^{d\times m}$ so that the data $A$ can be approximately represented by the reconstructed data $XX^\top A$. Its mathematical model can be expressed as
\begin{equation}\label{pca-c}
	\begin{aligned}
		\min_{X\in \mathbb{R}^{d\times m}}~& \|A-XX^\top A\|_\textrm{F}^2\\
		\rm{s.t.}~~~&X^\top X = I_m,
	\end{aligned}
\end{equation} 
where $m$ represents the transformed dimension. By expanding the objective function, the above problem can be equivalently transformed into
\begin{equation}\label{pca}
	\begin{aligned}
		\min_{X\in \mathbb{R}^{d\times m}}~& -\textrm{Tr}(X^\top A A^\top X)\\
		\rm{s.t.}~~~&X^\top X = I_m.
	\end{aligned}
\end{equation}

It is well demonstrated in \cite{wang2023simultaneous,xiu2022efficient} that $\ell_{2,0}$-norm is the most suitable for feature selection and can represent sparse structures more efficiently. In this regard, FSPCA selects discriminative features by adding $\ell_{2,0}$-norm in \eqref{pca}. The constrained SPCA model is given by
\begin{equation}\label{fspca-m}
	\begin{aligned}
		\min_{X\in \mathbb{R}^{d\times m}}~&- \textrm{Tr}(X^\top AA^\top X)\\
		\rm{s.t.}~~~&X^\top X = I_m,~\|X\|_{2,0}\leq r,
	\end{aligned}
\end{equation}
where $r>0$ indicates the number of non-zero rows corresponding to the features to be selected. Obviously, one can easily determine the number of selected features directly by adjusting the sparsity level $r$.

\section{The Proposed Method} \label{method}

This section first introduces our mathematical model, then presents the designed optimization algorithm in detail, and finally gives the rigorous convergence analysis.

\subsection{New Formulation}

\begin{figure}[t]
	\centering
	\includegraphics[width=0.9\textwidth]{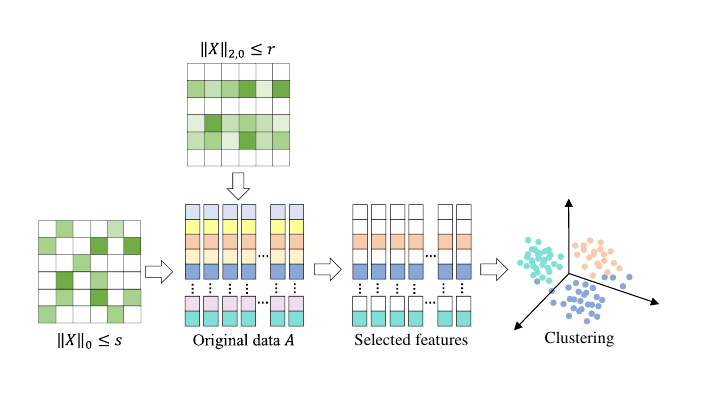}
	\vspace{-0.8cm}
	\caption{The flowchart of feature selection and clustering of our proposed DSCOFS, where $\|X\|_{2,0}$ captures the global structural sparsity, and $\|X\|_0$ captures the local element-wise sparsity.}
	\label{plot-1}
\end{figure}

Although $\ell_{2,0}$-norm can find the most significant structural features, the original local redundant elements in the data may affect the selection of structural features. In order to remove these local redundant elements, we further enforce the $\ell_0$-norm constraint to the transformation matrix $X$, and then obtain our DSCOFS model by considering the following double sparsity constrained optimization problem
\begin{equation}\label{dsco}
	\begin{aligned}
		\min_{X\in \mathbb{R}^{d\times m}}~&- \textrm{Tr}(X^\top A A^\top X)\\
		\rm{s.t.}~~~&X^\top X = I_m,~\|X\|_{2,0}\leq r,~\|X\|_0\leq s,  
	\end{aligned}
\end{equation}
where $s>0$ denotes the number of non-zero elements, which is used to characterize the element-wise sparsity. Compared with FSPCA in \eqref{fspca-m}, the advantages of our proposed DSCOFS  are
\begin{itemize}
\item $\|X\|_0\leq s$ helps in filtering out irregular noise and local individual sparsity, and then complements $\|X\|_{2,0}\leq r$ to more easily identify differential features. 
    \item Double sparsity provides more flexibility in feature selection because it allows variables $r$ and $s$ to be adjusted according to actual needs.
\end{itemize}
In addition, Figure \ref{plot-1} shows the flowchart of feature selection and clustering of our proposed DSCOFS.

\subsection{Optimization Algorithm}

One can easily observe that \eqref{dsco} has three nonconvex constraints, i.e., $X^\top X = I_m,~\|Y\|_0\leq s,~\|Z\|_{2,0}\leq r$, which are difficult to solve simultaneously. Motivated by \cite{chen2020alternating}, it is possible for us to update one variable while fixing other variables, which is called the alternating minimization algorithm  (AMA). To this end, \eqref{dsco} can be reformulated as 
\begin{equation}\label{dsco-1}
	\begin{aligned}
		\min_{X, Y, Z\in \mathbb{R}^{d\times m}}~&- \textrm{Tr}(X^\top A A^\top X)\\
		\rm{s.t.}~~~~~&X^\top X = I_m,~\|Y\|_0\leq s,~\|Z\|_{2,0}\leq r,\\
&X=Y,~X=Z,
	\end{aligned}
\end{equation}
where $Y$ and $Z$ are the introduced intermediate variables. For the sake of notation, let 
\begin{equation}
	\begin{aligned}
\mathcal{M}&=\{X\in\mathbb{R}^{d\times m} \mid X^\top X = I_m\},\\
\mathcal{S}&=\{Y\in\mathbb{R}^{d\times m} \mid \|Y\|_{0}\leq s\},\\
\mathcal{R}&=\{Z\in\mathbb{R}^{d\times m} \mid \|Z\|_{2,0}\leq r\}.
	\end{aligned}
\end{equation}
Then, \eqref{dsco-1} can be written as the penalty version
\begin{equation}\label{dsco-2}
	\begin{aligned}
		\min_{X, Y, Z\in \mathbb{R}^{d\times m}}~&- \textrm{Tr}(X^\top A A^\top X)+\mu_1 \|X-Y\|_\textrm{F}^2+\mu_2 \|X-Z\|_\textrm{F}^2\\
		\rm{s.t.}~~~~~&X\in\mathcal{M},~Y\in\mathcal{S},~Z\in\mathcal{R},\\
	\end{aligned}
\end{equation}
where $\mu_1, \mu_2>0$ are the penalty parameters. Overall, the iterative optimization framework of solving \eqref{dsco-2} (also \eqref{dsco}) can be presented in Algorithm \ref{am1}, where $\tau_1, \tau_2, \tau_3>0$ are proximal parameters to guarantee convergence.
\begin{algorithm}[t]
\caption{Optimization algorithm of solving  \eqref{dsco-2}} \label{am1}
\textbf{Input:} Data $A$, sparsity $s$, $r$, paraments $\mu_1$, $\mu_2$, $\tau_1$, $\tau_2$, $\tau_3$\\
\textbf{Initialize:} $(X^0, Y^0, Z^0)$\\
\textbf{While} not converged \textbf{do}	
		\begin{algorithmic}[1]
\STATE  Update  $X^{k+1}$ by
\begin{equation}\label{sub-x}
	\begin{aligned}
		\min_{X\in \mathbb{R}^{d\times m}}~& -\textrm{Tr}(X^\top A A^\top  X)+ \mu_1 \|X-Y^k\|_\textrm{F}^2+\mu_2 \|X-Z^k\|_\textrm{F}^2+\tau_1 \|X-X^k\|_\textrm{F}^2\\
		\textrm{s.t.}~~~&X\in\mathcal{M}
	\end{aligned}
\end{equation}
\STATE  Update $Y^{k+1}$ by
\begin{equation}\label{sub-y}
	\begin{aligned}
		\min_{Y\in \mathbb{R}^{d\times m}}~& \|X^{k+1}-Y\|_\textrm{F}^2+\tau_2 \|Y-Y^k\|_\textrm{F}^2 \\
		\textrm{s.t.}~~~&Y\in\mathcal{S}
	\end{aligned}
\end{equation}
\STATE  Update $Z^{k+1}$ by
\begin{equation}\label{sub-z}
	\begin{aligned}
		\min_{Z\in \mathbb{R}^{d\times m}}~& \|X^{k+1}-Z\|_\textrm{F}^2+\tau_3 \|Z-Z^k\|_\textrm{F}^2 \\
		\textrm{s.t.}~~~&Z\in \mathcal{R}
	\end{aligned}
\end{equation}
\STATE   Check convergence
		\end{algorithmic}
\textbf{End While}\\
\textbf{Output:} ($X^{k+1}, Y^{k+1}, Z^{k+1}$)
\end{algorithm}

In what follows, we focus on how to solve the resulting subproblems in detail.

\subsubsection{Update $X^{k+1}$}

Note that the constraint $X\in\mathcal{M}$ in \eqref{sub-x} is a Riemannian manifold, which is known as the orthogonal manifold. Although there are various optimization methods on orthogonal manifolds, most of them rely on generating points along the geodesics of $\mathcal{M}$, which makes the computation inefficient. To preserve the structure of the orthogonal constraint and avoid calculating the geodesics, we adopt an effective infeasibility optimization technique called the exact penalty function method \cite{xiao2022class}.

By introducing the merit function in \cite{gao2018new}, problem \eqref{sub-x}  can be transformed into the following equivalent problem containing a novel penalty function with a compact convex constraint $\mathcal{C}$, i.e.,
\begin{equation}\label{penc-x}
\begin{aligned}
\min_{X\in{\mathcal{C}}}~&h(X)=l(X)+g(X),
\end{aligned}
\end{equation}
where 
\begin{equation}
	\begin{aligned}
		l(X)=&-\textrm{Tr}(X^\top A A^\top  X)+ \mu_1 \|X-Y^k\|_\textrm{F}^2+\mu_2 \|X-Z^k\|_\textrm{F}^2+\tau_1 \|X-X^k\|_\textrm{F}^2
	\end{aligned}
\end{equation}
and 
\begin{equation}\label{penc-p}
g(X)=-\frac{1}{2}\langle \Lambda(X), X^\top X - I_m\rangle+\frac{\beta}{4}\|X^\top X - I_m\|_\textrm{F}^2,
\end{equation}
with the penalty parameter $\beta>0$ and 
\begin{equation}\label{penc-l}
\Lambda(X)=\frac{1}{2}(X^\top\nabla l(X)+ \nabla l(X)^\top X).
\end{equation}
The equivalence between problem \eqref{sub-x} and problem \eqref{penc-x} in the sense of global minimizers is established in  \cite[Theorem 3.2]{xiao2022class} by $\mathcal{C}$ chosen as a ball with radius $\rho$ in $\textrm{F}$-norm, i.e.,  $\mathcal{B}_\rho=\{X\in\mathbb {R}^{d\times m}\mid \|X\|_\textrm{F}\leq \rho\}$.
Since $\nabla l(X)$ is involved in $h(X)$, in order to avoid calculating the Hessian of $l(X)$ when calculating the gradient of $h(X)$, we consider an approximation of $\nabla h(X)$ as 

\begin{equation}\label{dx}
\begin{aligned}
D(X)= &-2AA^\top X +2\mu_1(X-Y^k)+2\mu_2(X-Z^k)\\&+2\tau_1(X-X^k)-X\Lambda(X)+\beta X(X^\top X - I_m).
\end{aligned}
\end{equation}
Then, \eqref{sub-x} can be solved by  the approximate gradient descent method for \eqref{penc-x}.  
One can derive the iterative formula of the gradient descent method given by
\begin{equation}\label{xk-1}
\begin{aligned}
\hat{X}^{k+1}=X^k -\eta_kD(X^k),
\end{aligned}
\end{equation}
where $D(X^k)$ is the approximate gradient of $h(X)$ at $X^k$ and $\eta_k>0$ is the Barzilai-Borwein step-size \cite{huang2021equipping}. 
\begin{algorithm}[t]
\caption{Exact penalty function method of solving \eqref{sub-x}} \label{pencf}
\textbf{Input:} Data $A$, paraments $\beta$, $\rho$, $\eta$\\
\textbf{Initialize:} $X^0 \leftarrow X^k$\\
\textbf{While} not converged \textbf{do}	
		\begin{algorithmic}[1]
\STATE  Compute $D(X^k)$ by \eqref{dx}
\STATE  Compute $\hat{X}^{k+1}$ by \eqref{xk-1}
\STATE  If $\|\hat{X}^{k+1}\|_\textrm{F}>\rho$ then\\
\begin{equation}\label{xk-11}
	\begin{aligned}
		X^{k+1}=\frac{\rho}{\|\hat{X}^{k+1}\|_\textrm{F}}\hat{X}^{k+1}
	\end{aligned}
\end{equation}
        else
\begin{equation}\label{xk-12}
	\begin{aligned}
		X^{k+1}=\hat{X}^{k+1}
	\end{aligned}
\end{equation}
        end if
\STATE   Check convergence
		\end{algorithmic}
\textbf{End While}\\
\textbf{Output:} $X^{k+1}$
\end{algorithm}

However, the above iteration cannot guarantee that $\hat{X}^{k+1}$ always satisfies the orthogonal manifold constraint $\mathcal{M}$. For this purpose, it chooses the  ball with radius $\rho>\sqrt{m}$. Denote $\mathcal{P}_ {\mathcal{B}_\rho}$ as the projection onto $\mathcal{B}_\rho$. Then, $X^{k+1}$ will be refined via
\begin{equation}\label{xk-2}
\begin{aligned}
X^{k+1}=\mathcal{P}_{\mathcal{B}_\rho}(\hat{X}^{k+1}).
\end{aligned}
\end{equation}

By simply shrinking $X^{k+1}$, it can be guaranteed that $X^{k+1}$ is within the range of $\mathcal{B}_\rho$. Of course, if $X^{k+1}$ is close enough to the orthogonal manifold constraint $\mathcal{M}$, one can abandon the calculation of $\mathcal{P}_{\mathcal{B}_\rho}(\cdot)$ and use \eqref{xk-1} to update $X^{k+1}$. The detailed iterative scheme for updating $X^{k+1}$ is provided in Algorithm \ref{pencf}.

\subsubsection{Update $Y^{k+1}$}
Once $X^{k+1}$ has been updated, it is time to solve \eqref{sub-y}. By direct calculation, problem \eqref{sub-y} can be written as 
\begin{equation}
\begin{aligned}
		\min_{Y\in \mathbb{R}^{d\times m}}~& \|\frac{X^{k+1}+\tau_2 Y^k}{1+\tau_2}-Y\|_\textrm{F}^2 \\
		\textrm{s.t.}~~~&Y\in\mathcal{S}.
	\end{aligned}
\end{equation}
 For ease of description, denote $W^{k+1}=\frac{X^{k+1}+\tau_2 Y^k}{1+\tau_2}$ and let $t^{k+1}_{s}>0$ be the $s$-th largest absolute value of all $W^{k+1}_{ij}$. By exploiting the hard thresholding operator \cite{blumensath2009iterative}, we obtain a closed-form solution given by
\begin{equation}\label{yk}
\begin{aligned}
Y^{k+1}_{ij}=\begin{cases}
   W^{k+1}_{ij},&~~|W^{k+1}_{ij}|\geq t^{k+1}_{s},\\
    0,&~~|W^{k+1}_{ij}|< t^{k+1}_{s},
\end {cases}
\end{aligned}
\end{equation}
where $|\cdot|$ denotes the absolute value. By setting different $s$, the solution $Y^{k+1}$ can be directly determined using the above truncation strategy.

\subsubsection{Update $Z^{k+1}$}
Similarly, \eqref{sub-z} can be reformulated as 
\begin{equation}
\begin{aligned}
		\min_{Z\in \mathbb{R}^{d\times m}}~& \|\frac{X^{k+1}+\tau_3 Z^k}{1+\tau_3 }-Z\|_\textrm{F}^2 \\
		\textrm{s.t.}~~~&Z\in\mathcal{R}.
	\end{aligned}
\end{equation}
Denote $V^{k+1}=\frac{X^{k+1}+\tau_3 Z^k}{1+\tau_3}$ and let $t^{k+1}_{r}>0$ represent the $r$-th largest element of all $\|V^{k+1}_{i:}\|_2$. Considering the structural sparsity of $\ell_{2,0}$-norm, \eqref{sub-z} has the following closed-form solution 
\begin{equation}\label{zk}
\begin{aligned}
Z^{k+1}_{i:}=\begin{cases}
    V^{k+1}_{i:},&~~\|V^{k+1}_{i:}\|_2\geq t^{k+1}_{r},\\
    0,&~~\|V^{k+1}_{i:}\|_2 < t^{k+1}_{r}.
\end {cases}
\end{aligned}
\end{equation}
In this regard, only the top-$r$ rows are retained, thus achieving the purpose of feature selection.

\begin{Remark}
From the perspective of optimization, the order of updating $Y^{k+1}$ and $Z^{k+1}$ in Algorithm \ref{am1} has no effect on the result. However, in the application of feature selection, we should first update $Y^{k+1}$ on $\mathcal{S}$ to remove redundancy, and then update $Z^{k+1}$ on $\mathcal{R}$ to select features. Otherwise, it is difficult to affect the results of feature selection. In addition, the optimization strategy developed in \cite{nie2022learning} is based on the rank of the covariance matrix, which is essentially different from ours.
\end{Remark}

\subsection{Convergence Analysis}

\label{convergence}
In this part, we will establish the global convergence for Algorithm \ref{am1} inspired by  \cite{bolte2014proximal}. For the sake of notation simplicity, we  denote the objective function of \eqref{dsco-2} as 
\begin{equation}\label{fxyz}
\begin{aligned}
f(X,Y,Z)=&-\textrm{Tr}(X^\top A A^\top X)+\mu_1 \|X-Y\|_\textrm{F}^2+\mu_2 \|X-Z\|_\textrm{F}^2.
\end{aligned}
\end{equation}
Here, $f(X,Y,Z)$ is continuously differentiable and its gradient is given by
\begin{equation}\label{fgra}
\begin{aligned}
\nabla f(X, Y, Z)&=\left[\begin{array}{c}
\frac{\partial}{\partial X}f(X, Y, Z)\\
\frac{\partial}{\partial Y}f(X, Y, Z)\\
\frac{\partial}{\partial Z}f(X, Y, Z)
\end{array}\right]\\
&=\left[\begin{array}{c}
-2(AA^\top X -\mu_1(X-Y)-\mu_2(X-Z))\\
2\mu_1(Y-X)\\
2\mu_2(Z-X)
\end{array}\right].\end{aligned}
\end{equation}

Let $ \textrm{N}_{\mathcal{M}\times\mathcal{S}\times\mathcal{R}} (X, Y, Z)$ be the normal cones of $\mathcal{M}\times\mathcal{S}\times\mathcal{R}$ at $(X, Y, Z)$ and denote 
\begin{equation}
\begin{aligned}\label{SCA1}
\lambda_0 &=\sup_{X\in \mathcal{B}_\rho}\max\{1,\|\nabla l(X)\|_\textrm{F}\},\\
 \lambda_1 &=\sup_{X\in \mathcal{B}_\rho}\max\{1,\|\Lambda(X)\|_\textrm{F}\}, \\
 \lambda_2 &=\sup_{X_1,X_2\in \mathcal{B}_\rho
  }\max\{1,\frac{\|\Lambda(X_1)-\Lambda(X_2)\|}{\|X_1-X_2\|}\},
\end{aligned}
\end{equation}
where $\|\cdot\|$ represents the spectral norm, i.e., the largest singular value.
%

\begin{Theorem}\label{Theor1}
Suppose that $\beta\geq \max\{2(\lambda_0+\lambda_1),2m\lambda_2\}$. Let $\{(X^k, Y^k, Z^k)\}$ be the sequence generated by Algorithm \ref{am1} for solving problem \eqref{dsco-2}. Then the following properties hold: 
\begin{itemize}
  \item [(a)] $\{f(X^k, Y^k, Z^k)\}$ is strictly nonincreasing;
  \item [(b)] The sequence $\{(X^k, Y^k, Z^k)\}$ is bounded;
  \item [(c)] $\lim_{k\rightarrow\infty}\|(X^{k+1}, Y^{k+1}, Z^{k+1})-(X^k, Y^k, Z^k)\|_\textrm{F}=0$;
  \item [(d)] Any accumulation point $(X^*, Y^*, Z^*)$ of the sequence $\{(X^k, Y^k, Z^k)\}$ is a stationary point of \eqref{dsco-2} in the sense that 
  $$0\in\nabla f(X^*, Y^*, Z^*))+\textrm{N}_{\mathcal{M}\times\mathcal{S}\times\mathcal{R}} (X^*, Y^*, Z^*)). $$
\end{itemize}
\end{Theorem} 
The proofs  are provided in the Appendix.

\begin{Remark}
The above theorem shows the generated sequence by Algorithm \ref{am1} is strictly non-increasing, which shares a similar result in \cite{li2021sparse}. Fortunately, with the help of Kurdyka-Lojasiewicz (KL) properties, we further establish the global convergence without any further assumption, which is different from \cite{nie2022learning}.
\end{Remark}

\section{Numerical Experiments} \label{expriments}

This section validates the effectiveness and superiority of our proposed DSCOFS by comparing with benchmark methods, including LapScore \cite{he2005laplacian}, UDFS \cite{yang2011l}, SOGFS \cite{nie2019structured}, RNE \cite{liu2020robust}, FSPCA \cite{nie2022learning}, SPCAFS \cite{li2021sparse}, and SPCA-PSD \cite{zheng2023fast}. In particular, LapScore, UDFS, SOGFS, and RNE are implemented using Auto-UFSTool\footnote{https://github.com/farhadabedinzadeh/AutoUFSTool},  FSPCA\footnote{https://github.com/tianlai09/FSPCA}, SPCAFS\footnote{https://github.com/quiter2005/algorithm}, and SPCA-PSD\footnote{https://github.com/zjj20212035/SPCA-PSD} are implemented using the codes downloaded from their provided links.

\begin{table}[t]
\caption{The dataset information.}\label{data}
    \vspace{0.2cm}
\centering
\begin{tabular}{|c|c|c|c|c|c|}
    \hline
    Type&Datasets & Features &Samples &Classes \\
    \hline\hline
    \multirow{3}*{Synthetic}&
     2Spiral&9&1000&2\\
    &Banana&9&1000&2\\
    &Dartboard&9&1000&4\\
    \hline
    \multirow{8}*{Real-world}&  
     COIL20 & 1024 & 1440 & 20 \\
    &USPS & 256 & 1000 & 10 \\
    &lung\_discrete & 325 & 73 & 7 \\
    &GLIOMA & 4434 & 50 & 4 \\
    &UMIST\tablefootnote{https://github.com/ycwang-libra/LpCNMF\_github/tree/main/datasets} & 644 & 575 & 20 \\
    &warpPIE10P & 2420 & 210 & 10 \\
    &Isolet & 617 & 1560 & 26 \\
    &MSTAR\_SOC\_CNN\tablefootnote{https://github.com/zjj20212035/SPCA-PSD} & 1024 & 2425 & 10 \\
    \hline
\end{tabular}
\end{table}

Section \ref{exp-1} gives the experimental setup, Section \ref{exp-2} shows the quantitative results of feature selection, Section \ref{exp-3} analyzes the differences of double sparsity, and Section \ref{exp-4} provides more discussion.

\subsection{Experimental Setup}\label{exp-1}

\subsubsection{Dataset Description}

In the experiments, three synthetic datasets\footnote{https://github.com/milaan9/Clustering-Datasets/} and eight real-world datasets\footnote{https://jundongl.github.io/scikit-feature/datasets.html} are considered. The real-world datasets contain three face image datasets, i.e., COIL20, warpPIE10P, UMIST, one handwritten image dataset, i.e., USPS, two biological datasets, i.e., GLIOMA, lung\_discrete, one spoken letter recognition dataset, i.e., Isolet, and one deep learning dataset, i.e., MSTAR\_SOC\_CNN. Note that the lung\_discrete dataset is discrete, while the rest are continuous. For details of these datasets, see Table \ref{data}.

\subsubsection{Parameter Setting}

To be fair, it specifies the parameters of all comparison methods. 
For LapScore, the heatkernel connection is used to construct the similarity matrix and set $t=1$. 
For LapScore, UDFS, SOGFS, and RNE, the number of neighbors is set to 5 for constructing the weighted matrix. 
For UDFS and SOGFS, the number of clusters is set to the number of data categories. 
For SPCAFS, $p=0.5$ is chosen for the $\ell_{2,p}$-norm regularization. 
For SOGFS, SPCAFS, and our proposed DSCOFS, the projection dimension is fixed to the number of data categories.

According to \cite{li2021sparse}, the regularization parameters of all comparison methods are tuned from the candidate set
$\{10^{-6}, 10^{-4}, 10^{-2}, 10^0, 10^2, 10^4, 10^ 6\}$ using the grid search strategy. For all datasets, the number of selected features is chosen from $\{10, 20, 30, 40,$\\
$50, 60, 70, 80, 90, 100\}$.  
For our proposed DSCOFS, the element-wise sparsity parameter is set to $s=\alpha d m$, where $\alpha$ is the sparsity percentage, indicating the percentage of elements retained, and $d m$ represents the total number of elements in the transformation matrix. In this study, $\alpha$ is selected from  $\{0.1, 0.2, 0.3, 0.4, 0.5, 0.6, 0.7, 0.8, 0.9\}$.

\subsubsection{Initialization and Stopping Criteria}

In the case of DSOCFS, a random orthogonal matrix is employed as the initial solution. To reduce the impact of randomness, we perform 10 times and select the orthogonal matrix that minimizes $-\textrm{Tr}(X^\top A A^\top X)$ as our initial solution.

For Algorithm \ref{am1}, it is terminated if 
\begin{equation}\label{stop-dsco}
\begin{aligned}
\frac{|f(X^{k+1},Y^{k+1},Z^{k+1})-f(X^{k},Y^{k},Z^{k})|}{1+|f(X^{k},Y^{k},Z^{k})|}\leq 10^{-3}
\end{aligned}
\end{equation}
holds or the number of iterations reaches 100. In addition, the stopping criteria of Algorithm \ref{pencf} is consistent with \cite{xiao2022class}.

\subsubsection{Evaluation Metrics}

To evaluate the feature selection performance, the required features are initially obtained by all comparison methods. Based on the selected feature subspace, $K$-means clustering is conducted to achieve pseudo labels, and the Kuhn-Munkres algorithm is employed to identify the optimal correspondence between the pseudo and real labels. Then two popular metrics, i.e., clustering accuracy (ACC) and normalized mutual information (NMI), are selected to measure the clustering performance. Note that different initial values will result in different clustering performance, we perform $K$-means clustering 50 times and record the quantitative results with the best parameters.

ACC reflects the accuracy of clustering, which is defined as
\begin{equation}\label{acc}
\begin{aligned}
\textrm{ACC}=\frac{1}{n}\sum_{i=1}^n \delta(\phi(i), \varphi(i)),
\end{aligned}
\end{equation}
where $\phi(i)$ and $\varphi(i)$ represent the $i$-th pseudo label matched by the Kuhn-Munkres algorithm and the real label, respectively. Here, $\delta(\phi(i), \varphi(i))=1$ if $\phi(i)=\varphi(i)$, otherwise $\delta(\phi(i), \varphi(i))=0$.

NMI illustrates the similarity between the clustering results and the true results, which is defined as
\begin{equation}\label{nmi}
\begin{aligned}
\textrm{NMI}=\frac{I(P,Q)}{\sqrt{H(P)H(Q)}},
\end{aligned}
\end{equation}
where $P, Q$ are the pseudo label and the real label, respectively, $I(P, Q)$ represents their mutual information, and $H(P), H(Q)$ denote their entropy.
The larger the ACC and NMI, the better the clustering results.

\subsection{Feature Selection Experiments}\label{exp-2}

\begin{figure}[t]
	\centering
	\hspace{-2.5mm}
	\subfigcapskip=-5pt
	\subfigure[2Spiral]{
		\centering
		\label{a}
		\includegraphics[width=0.18\textwidth]{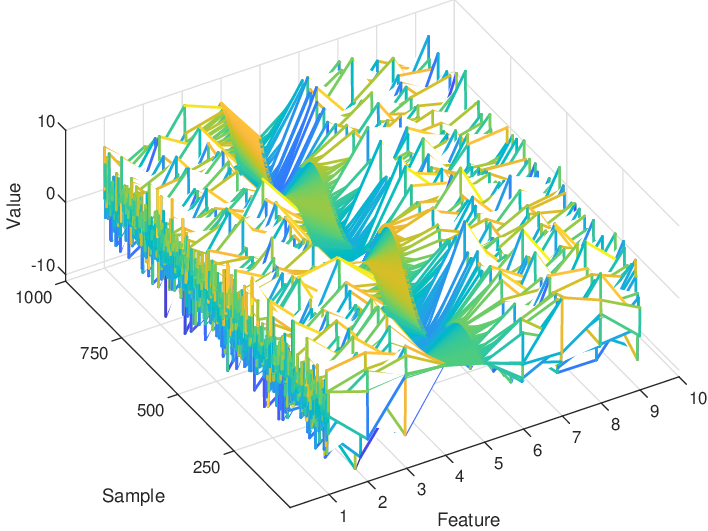}
	}
	\hspace{-1.5mm}
	\subfigcapskip=-5pt
	\subfigure[FSPCA]{		
		\label{b}
		\centering
		\includegraphics[width=0.18\textwidth]{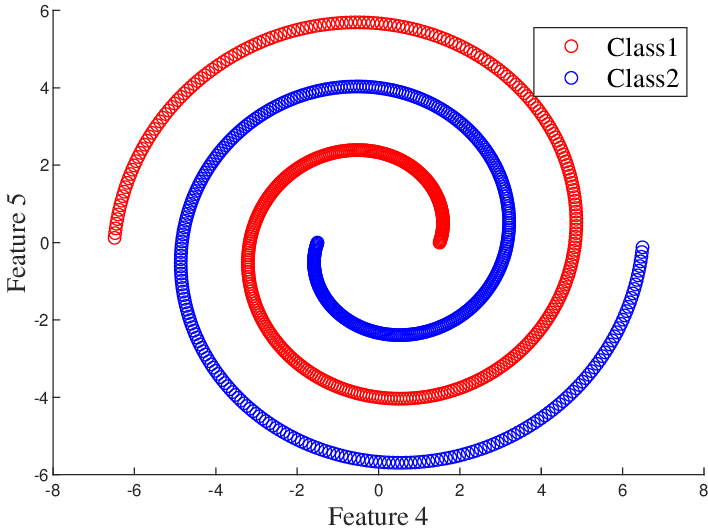}				
	}
	\hspace{-1.5mm}
	\subfigcapskip=-5pt
	\subfigure[SPCAFS]{
		\label{c}
		\centering
		\includegraphics[width=0.18\textwidth]{figure/2sprialf45}
	}
	\hspace{-1.5mm}
	\subfigcapskip=-5pt
	\subfigure[SPCA-PSD]{
		\label{d}
		\centering
		\includegraphics[width=0.18\textwidth]{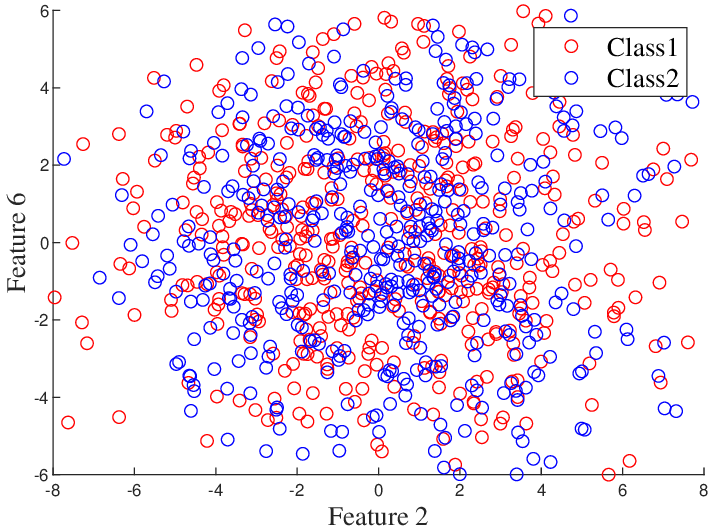}
	}
	\hspace{-1.5mm}
	\subfigcapskip=-5pt
	\subfigure[DSCOFS]{
		\label{d}
		\centering
		\includegraphics[width=0.18\textwidth]{figure/2sprialf45}
	}
	\hspace{-2.5mm}
	
	\hspace{-2.5mm}
	\subfigcapskip=-5pt
	\subfigure[Banana]{
		\centering
		\label{a}
		\includegraphics[width=0.18\textwidth]{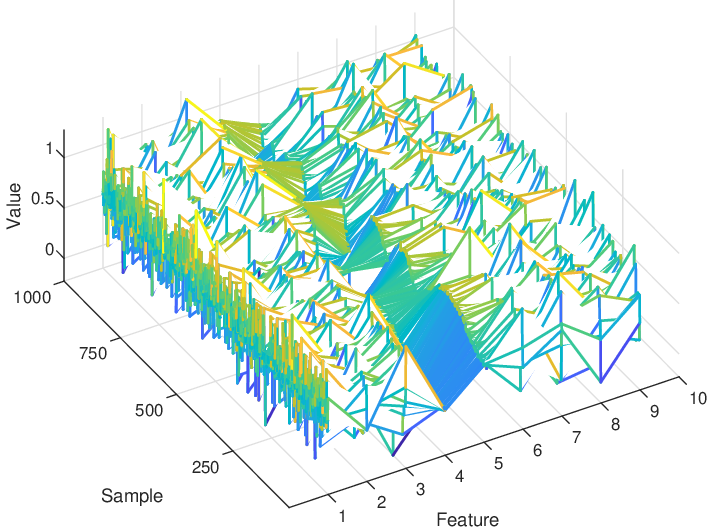}
	}
	\hspace{-1.5mm}
	\subfigcapskip=-5pt
	\subfigure[FSPCA]{		
		\label{b}
		\centering
		\includegraphics[width=0.18\textwidth]{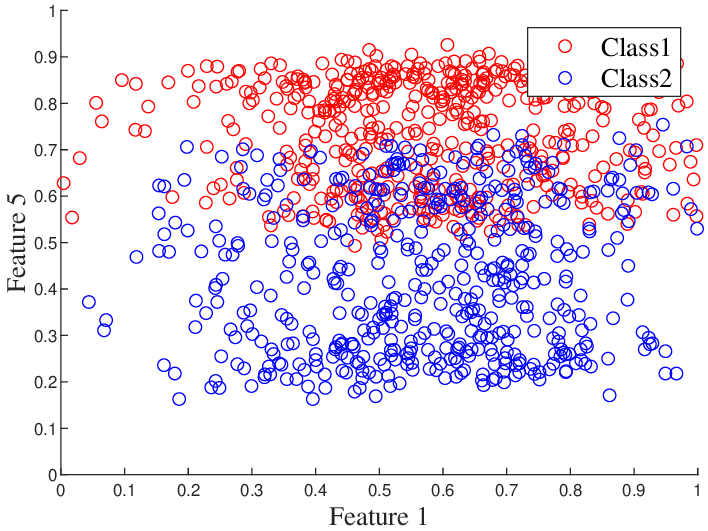}				
	}
	\hspace{-1.5mm}
	\subfigcapskip=-5pt
	\subfigure[SPCAFS]{
		\label{c}
		\centering
		\includegraphics[width=0.18\textwidth]{figure/banana15}
	}
	\hspace{-1.5mm}
	\subfigcapskip=-5pt
	\subfigure[SPCA-PSD]{
		\label{d}
		\centering
		\includegraphics[width=0.18\textwidth]{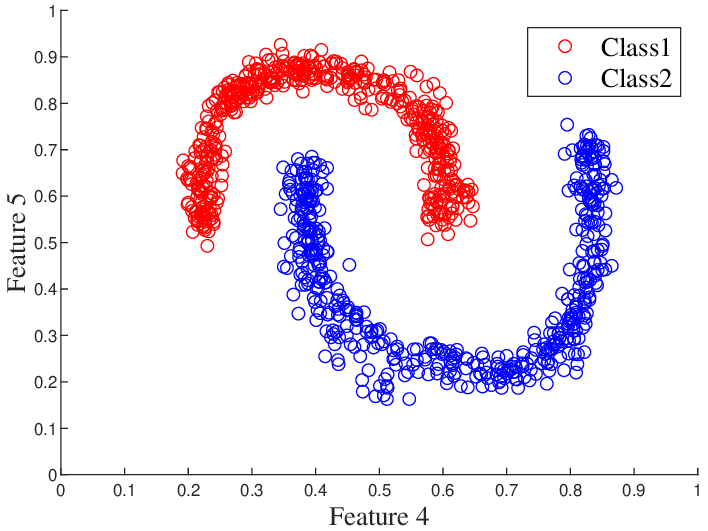}
	}
	\hspace{-1.5mm}
	\subfigcapskip=-5pt
	\subfigure[DSCOFS]{
		\label{d}
		\centering
		\includegraphics[width=0.18\textwidth]{figure/banana45}
	}
	\hspace{-2.5mm}
	
	\hspace{-2.5mm}
	\subfigcapskip=-5pt
	\subfigure[Dartboard]{
		\centering
		\label{a}
		\includegraphics[width=0.18\textwidth]{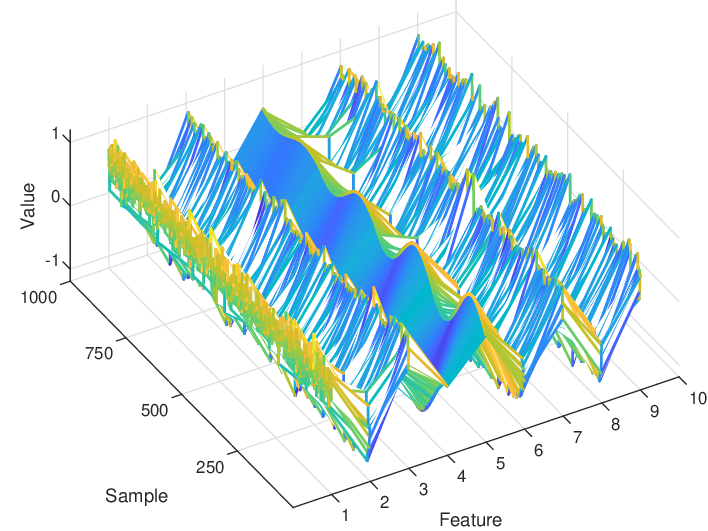}
	}
	\hspace{-1.5mm}
	\subfigcapskip=-5pt
	\subfigure[FSPCA]{		
		\label{b}
		\centering
		\includegraphics[width=0.18\textwidth]{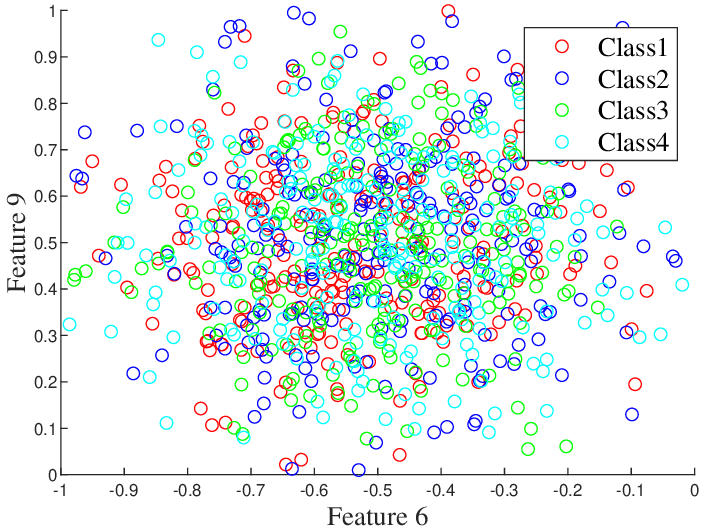}				
	}
	\hspace{-1.5mm}
	\subfigcapskip=-5pt
	\subfigure[SPCAFS]{
		\label{c}
		\centering
		\includegraphics[width=0.18\textwidth]{figure/dartboard69}
	}
	\hspace{-1.5mm}
	\subfigcapskip=-5pt
	\subfigure[SPCA-PSD]{
		\label{d}
		\centering
		\includegraphics[width=0.18\textwidth]{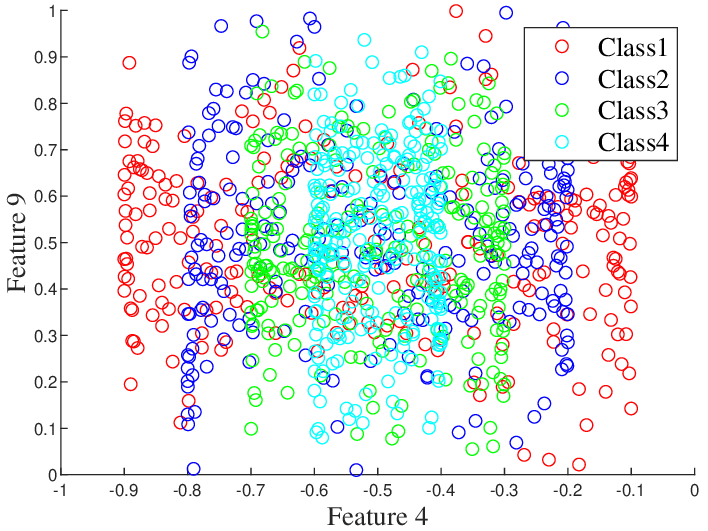}
	}
	\hspace{-1.5mm}
	\subfigcapskip=-5pt
	\subfigure[DSCOFS]{
		\label{d}
		\centering
		\includegraphics[width=0.18\textwidth]{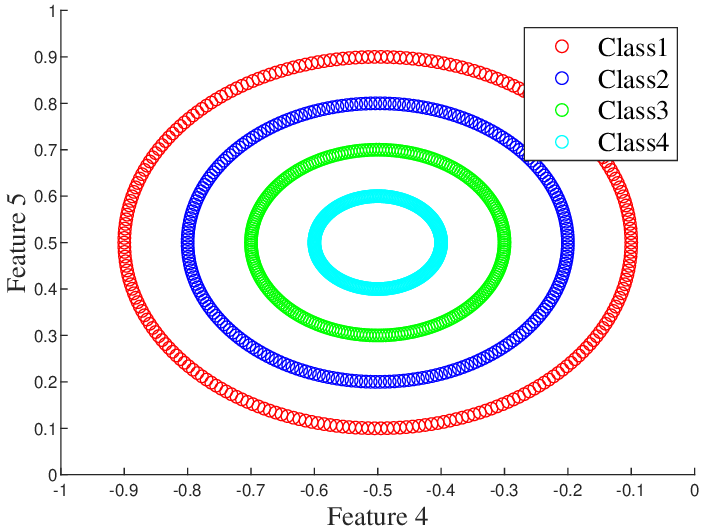}
	}
	\hspace{-2.5mm}
	\vspace{-0.2cm}
	\caption{The original distribution and feature selection results on the synthetic datasets. (a), (f) and (k) are the selected synthetic datasets; (b)-(e) are the feature selection results on 2Spiral; (g)-(j) are the feature selection results on Banana; (l)-(o) are the feature selection results on Dartboard. The top two features are used to show the results of feature selection and the correct features are feature 4 and feature 5.}
	\centering
	\label{plot-Synthetic}
\end{figure}

\begin{figure}[t]
	\centering
	\hspace{-4.5mm}
	\subfigcapskip=-5pt
	\subfigure[COIL20]{
		\centering
		\label{a}
		\includegraphics[width=0.23\textwidth]{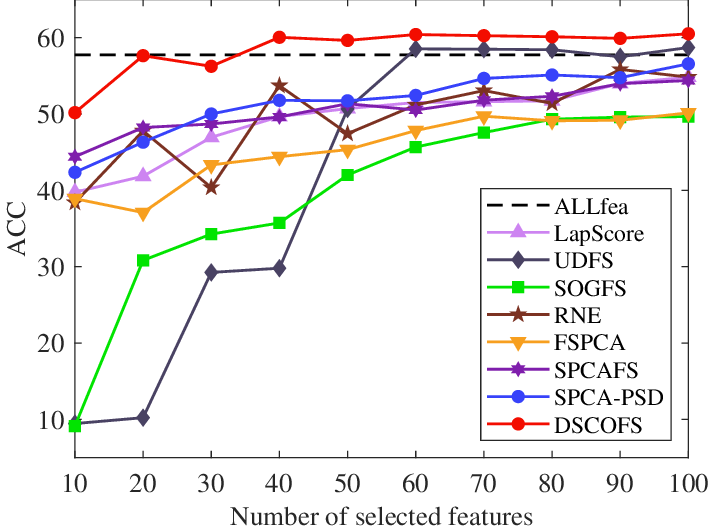}
	}   
	\hspace{-1.5mm} 
	\subfigcapskip=-5pt
	\subfigure[USPS]{		
		\label{b}
		\centering
		\includegraphics[width=0.23\textwidth]{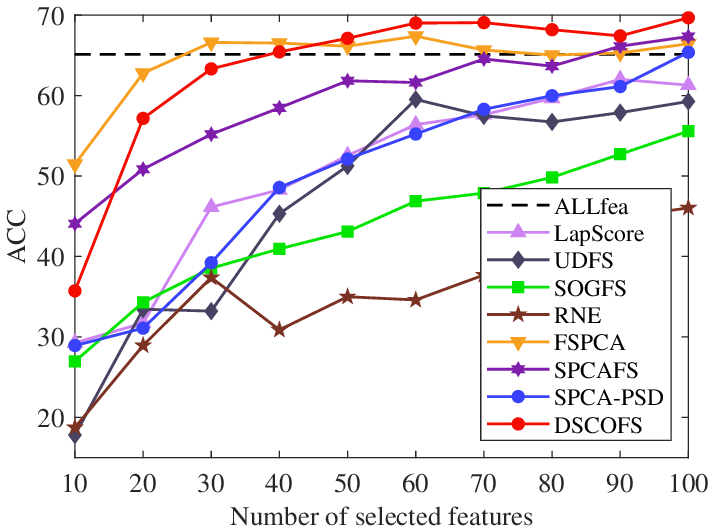}				
	}
	\hspace{-1.5mm} 
	\subfigcapskip=-5pt
	\subfigure[lung\_discrete]{
		\label{c}
		\centering
		\includegraphics[width=0.23\textwidth]{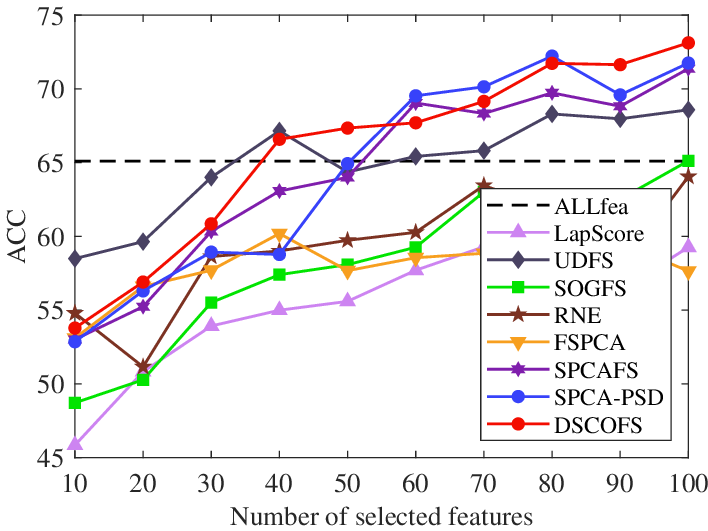}
	}
	\hspace{-1.5mm} 
	\subfigcapskip=-5pt
	\subfigure[GLIOMA]{
		\label{d}
		\centering
		\includegraphics[width=0.23\textwidth]{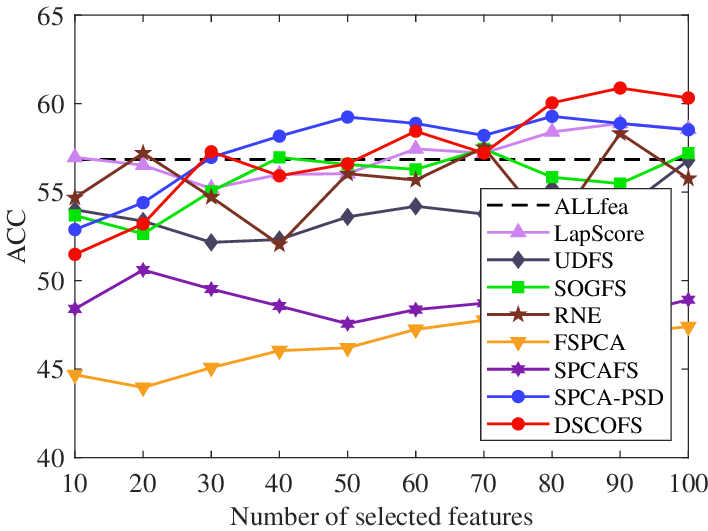}
	}
	\hspace{-4.5mm}
	
	\hspace{-4.5mm}
	\subfigcapskip=-5pt
	\subfigure[UMIST]{
		\label{e}
		\centering
		\includegraphics[width=0.23\textwidth]{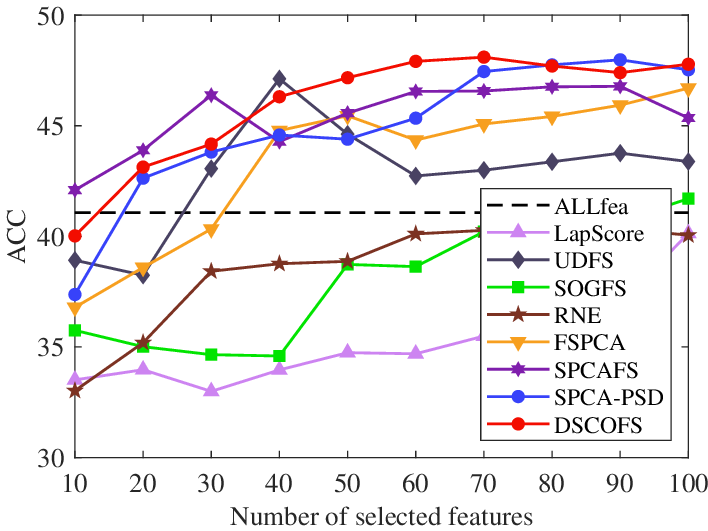}
	}
	\hspace{-1.5mm}  
	\subfigcapskip=-5pt
	\subfigure[warpPIE10P]{
		\label{f}
		\centering
		\includegraphics[width=0.23\textwidth]{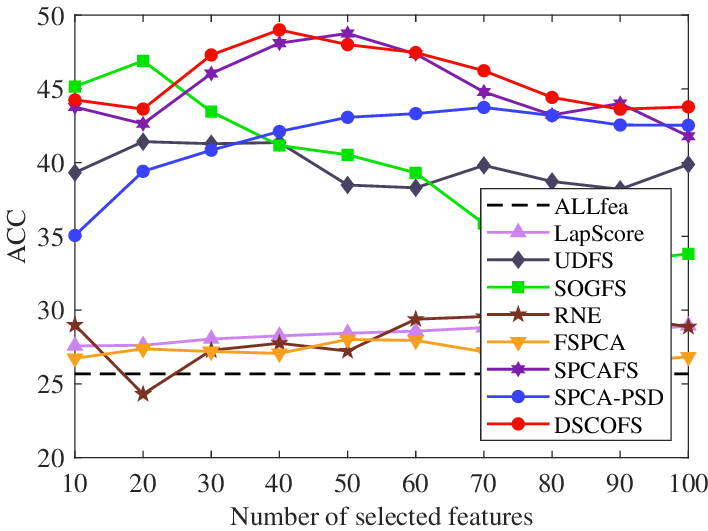}
	}
	\hspace{-1.5mm} 
	\subfigcapskip=-5pt
	\subfigure[Isolet]{
		\label{g}
		\centering
		\includegraphics[width=0.23\textwidth]{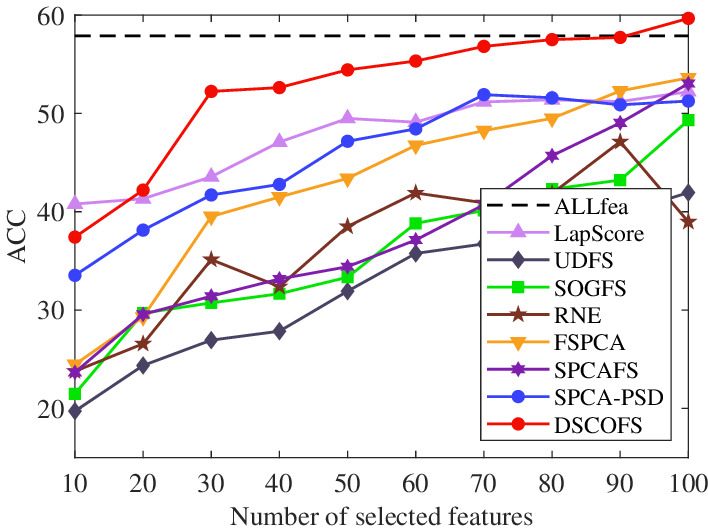}
	}
	\hspace{-1.5mm} 
	\subfigcapskip=-5pt
	\subfigure[MSTAR\_SOC\_CNN]{
		\label{g}
		\centering
		\includegraphics[width=0.23\textwidth]{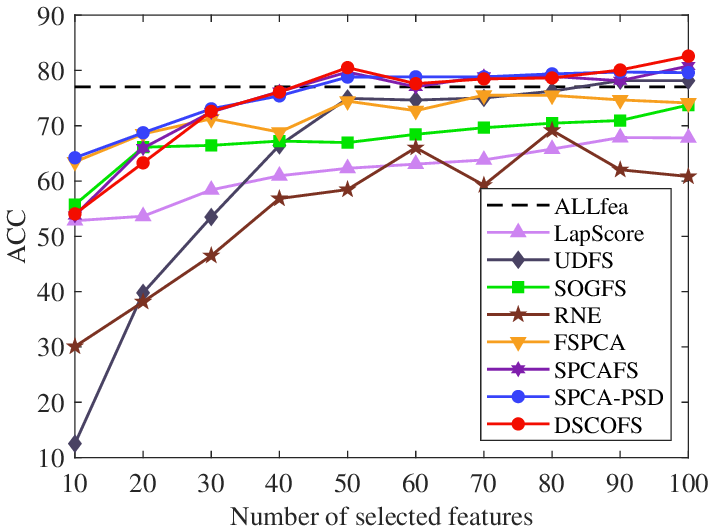}
	}
	\hspace{-4.5mm}
	\vspace{-0.2cm}
	\caption{The ACC (\%) curves  of compared methods on eight real-world datasets. All methods select features based on the given parameters and only the mean results of ACC (\%) are plotted.}
	\centering
	\label{plot-acc}
\end{figure}

\begin{table}[t]
	\caption{The ACC (mean\% $\pm$ std\%)  results  of compared methods on eight real-world datasets. The best and second-best results, except ALLfea, are marked in \textcolor[RGB]{255,0,0}{red} and \textcolor[RGB]{0,0,255}{blue}, respectively.}
	\label{Tacc}
	\vspace{0.2cm}
	\centering
	\scalebox{0.55}{
		\begin{tabular}{|c|c|c|c|c|c|c|c|c|c|c|}
			\hline
			Datasets& ALLfea & LapScore & UDFS & SOGFS & RNE & FSPCA & SPCAFS & SPCA-PSD & DSCOFS  \\
			\hline\hline
			\multirow{2}*{COIL20}&
			\multirow{2}*{57.74$\pm$4.93}
			& 54.82$\pm$3.91 & \textcolor[RGB]{0,0,255}{58.71$\pm$3.47} & 49.66$\pm$4.81 & 55.84$\pm$4.41 & 50.15$\pm$4.70 & 54.39$\pm$3.67 & 56.57$\pm$3.05 & \textcolor[RGB]{255,0,0}{60.51$\pm$4.63}\\
			&& (100) & \textcolor[RGB]{0,0,255}{(100)} & (100) & (90) & (100) & (100) & 100) & \textcolor[RGB]{255,0,0}{(100)}\\  
			\multirow{2}*{USPS}&
			\multirow{2}*{65.12$\pm$4.95} 
			& 62.02$\pm$4.09 & 59.52$\pm$2.97 & 55.58$\pm$3.07 & 46.04$\pm$2.69 & \textcolor[RGB]{0,0,255}{67.38$\pm$4.36} & 67.34$\pm$4.49 & 65.38$\pm$4.26 & \textcolor[RGB]{255,0,0}{69.67$\pm$4.97}\\
			&& (90) & (60) & (100) & (100) & \textcolor[RGB]{0,0,255}{(60)} & (100) & (100) & \textcolor[RGB]{255,0,0}{(100)}\\  
			\multirow{2}*{lung\_discrete}&
			\multirow{2}*{65.10$\pm$6.44} 
			& 59.29$\pm$6.33 & 68.58$\pm$6.99 & 65.12$\pm$6.89 & 64.05$\pm$6.65 & 60.19$\pm$6.55 & 71.37$\pm$7.68 & \textcolor[RGB]{0,0,255}{72.22$\pm$8.02} & \textcolor[RGB]{255,0,0}{73.12$\pm$8.48}\\
			&& (70) & (100) & (100) & (100) & (40) & (100) & \textcolor[RGB]{0,0,255}{(80)} & \textcolor[RGB]{255,0,0}{(100)}\\  
			\multirow{2}*{GLIOMA}&
			\multirow{2}*{56.84$\pm$5.24}
			&  58.88$\pm$3.96 & 56.80$\pm$4.85 & 57.44$\pm$6.16 & 58.32$\pm$7.31 & 47.92$\pm$4.61 & 50.60$\pm$5.02 & \textcolor[RGB]{0,0,255}{59.28$\pm$5.01} & \textcolor[RGB]{255,0,0}{60.88$\pm$6.31}\\
			&&  (90) & (100) & (70) & (90) & (80) & (20) & \textcolor[RGB]{0,0,255}{(90)} & \textcolor[RGB]{255,0,0}{(80)}\\ 
			\multirow{2}*{UMIST}&
			\multirow{2}*{41.07$\pm$2.38}           
			& 40.13$\pm$2.79 & 47.12$\pm$2.49 & 41.70$\pm$3.17 & 40.35$\pm$2.26 & 46.70$\pm$2.29 & 46.78$\pm$2.51 & \textcolor[RGB]{0,0,255}{47.98$\pm$2.91} & \textcolor[RGB]{255,0,0}{48.10$\pm$3.01}\\
			&& (100) & (40) & (100) & (90) & (100) & (90) & \textcolor[RGB]{0,0,255}{(90)} & \textcolor[RGB]{255,0,0}{(70)}\\  
			\multirow{2}*{warpPIE10P}&
			\multirow{2}*{25.67$\pm$1.90}
			& 28.94$\pm$1.66 & 41.42$\pm$3.18 & 46.90$\pm$3.89 & 29.57$\pm$2.96 & 28.01$\pm$2.27 & \textcolor[RGB]{0,0,255}{48.76$\pm$3.86} & 43.74$\pm$3.91 & \textcolor[RGB]{255,0,0}{49.00$\pm$3.88}\\
			&& (100) & (20) & (20) & (90) & (50) & \textcolor[RGB]{0,0,255}{(50)} & (70) & \textcolor[RGB]{255,0,0}{(40)}\\  
			\multirow{2}*{Isolet}&
			\multirow{2}*{57.89$\pm$3.82}
			& 52.21$\pm$2.76 & 41.95$\pm$2.07 & 49.31$\pm$2.32 & 47.12$\pm$2.06 & \textcolor[RGB]{0,0,255}{53.62$\pm$2.36} & 53.04$\pm$2.33 & 51.91$\pm$2.15 & \textcolor[RGB]{255,0,0}{59.67$\pm$3.46}\\
			&& (100)& (100) & (100) & (90) & \textcolor[RGB]{0,0,255}{(100)} & (100) & (70) & \textcolor[RGB]{255,0,0}{(100)}\\    
			\multirow{2}*{MSTAR\_SOC\_CNN}&
			\multirow{2}*{77.04$\pm$7.98}
			& 67.87$\pm$3.49 & 78.15$\pm$5.80 & 73.74$\pm$5.89 & 69.16$\pm$6.03 & 75.52$\pm$6.22 & \textcolor[RGB]{0,0,255}{80.80$\pm$5.95} & 79.70$\pm$6.43 & \textcolor[RGB]{255,0,0}{82.59$\pm$7.41}\\
			&& (90) & (90) & (100) & (100) & (70) & \textcolor[RGB]{0,0,255}{(100)} & (90) & \textcolor[RGB]{255,0,0}{(100)}\\  
			\hline
			Average& 55.81$\pm$4.71 & 53.02$\pm$3.62 & 56.53$\pm$4.04 & 54.93$\pm$4.53 & 51.31$\pm$4.30 & 53.69$\pm$4.47 & 59.14$\pm$4.44 & \textcolor[RGB]{0,0,255}{59.60$\pm$4.47} & \textcolor[RGB]{255,0,0}{62.94$\pm$5.27}\\
			\hline
		\end{tabular}
	}
\end{table}

\subsubsection{Comparison on Synthetic Datasets}

In this experiment, the synthetic datasets originally have only two features. The mean and variance of the original features are used to generate the remaining seven Gaussian noise features, and the original features are placed at the positions of features 4 and 5. The Banana dataset consists of 500 samples randomly selected from each of the two classes of the original samples. The SPCA-based feature selection methods are chosen, including FSPCA, SPCAFS, SPCA-PSD, and our proposed DSCOFS. It should be noted that FSPCA considers $\ell_{2,0}$-norm, SPCAFS considers $\ell_{2,p}$-norm $(0<p\leq 1)$, and SPCA-PSD considers $\ell_{2,1}$-norm, while our proposed DSCOFS considers double sparsity, i.e., $\ell_{2,0}$-norm plus $\ell_0$-norm.  

Figure \ref{plot-Synthetic} shows the feature selection results on the above three synthetic datasets, where the first column visualizes the original datasets and the last three columns visualize the results obtained by different methods. The top two features are used to show the results and feature 4 and feature 5 are correct. It can be seen that our proposed DSCOFS selects the correct features on all three datasets. For FSPCA and SPCAFS, the correct features are only selected on 2Spiral, while SPCA-PSD only selects the correct features on Banana. However, all methods except our proposed DSCOFS cannot select the correct features on  Dartboard. Obviously, different sparsity will lead to different results in feature selection. In contrast, our proposed DSCOFS takes both structural sparsity and element-wise sparsity into account, providing more possibilities for feature selection.

\subsubsection{Comparison on Real-world Datasets}
For eight real-world datasets,  Figure \ref{plot-acc} shows the mean curves of ACC under different feature numbers, where ALLfea represents all features used for clustering, and the results are considered as the baseline for analysis. In addition, Table \ref{Tacc} presents the mean and standard deviation of the best ACC within 100 features and the corresponding number of selected features. 
It can be seen from Figure \ref{plot-acc} that our proposed DSCOFS is the only method that outperforms the baseline on all datasets. Compared with other comparison methods, the ACC curves of our proposed DSCOFS show excellent performance, achieving the best results on all datasets, especially on COIL20 and Isolet. From Table \ref{Tacc}, it can be concluded that the average ACC of all datasets has improved by 3.34\%. For COIL20, UDFS outperforms FSPCA, SPCAFS, and SPCA-PSD, but not as well as our proposed DSCOFS, which improves by 1.8\%. For Isolet, our proposed DSCOFS is 6.05\%, 6.63\%, and 7.76\% better than FSPCA, SPCAFS, and SPCA-PSD, respectively.

\begin{figure}[t]
	\centering
	\hspace{-4.5mm}
	\subfigcapskip=-5pt
	\subfigure[COIL20]{
		\centering
		\label{a}
		\includegraphics[width=0.23\textwidth]{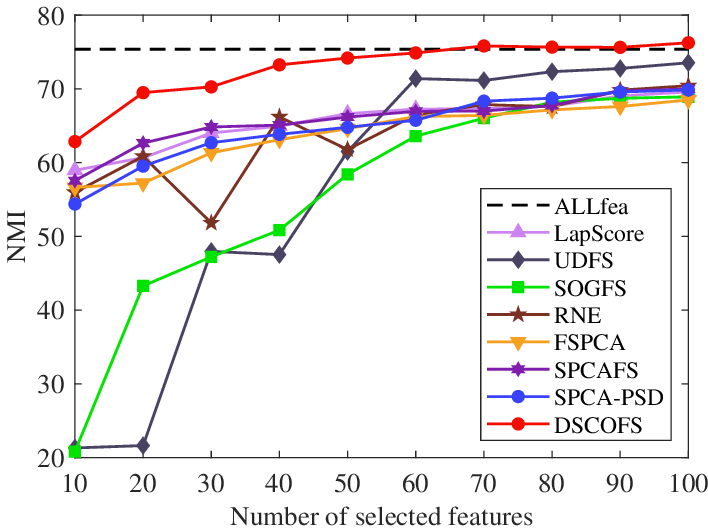}
	}
	\hspace{-1.5mm} 
	\subfigcapskip=-5pt
	\subfigure[USPS]{		
		\label{b}
		\centering
		\includegraphics[width=0.23\textwidth]{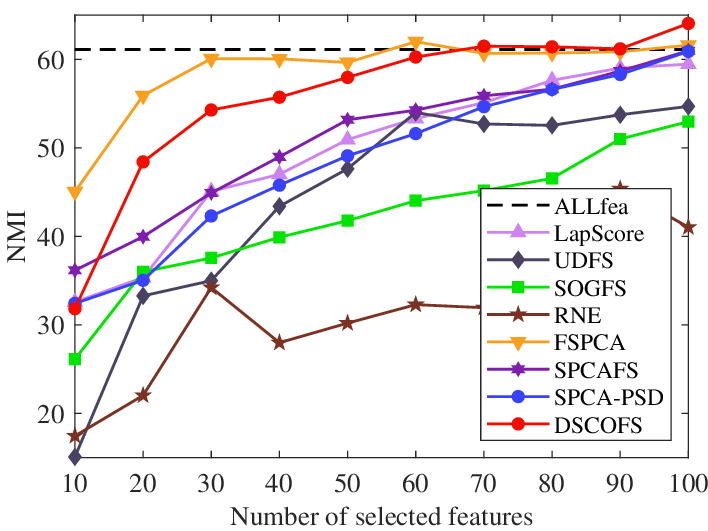}				
	}
	\hspace{-1.5mm} 
	\subfigcapskip=-5pt
	\subfigure[lung\_discrete]{
		\label{c}
		\centering
		\includegraphics[width=0.23\textwidth]{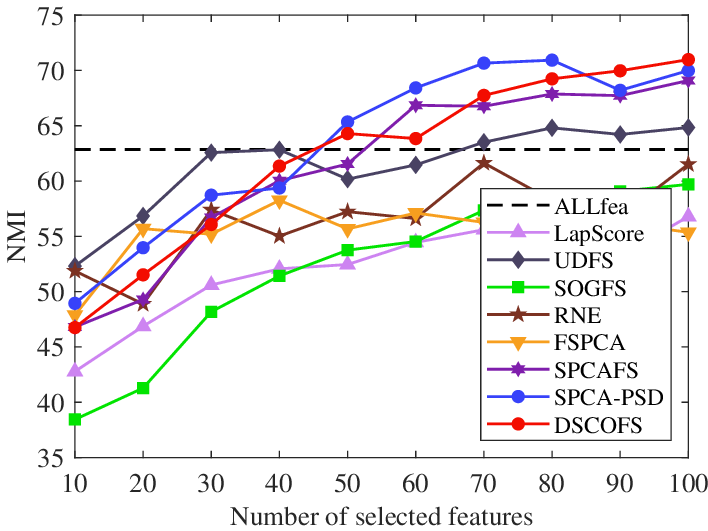}
	}
	\hspace{-1.5mm} 
	\subfigcapskip=-5pt
	\subfigure[GLIOMA]{
		\label{d}
		\centering
		\includegraphics[width=0.23\textwidth]{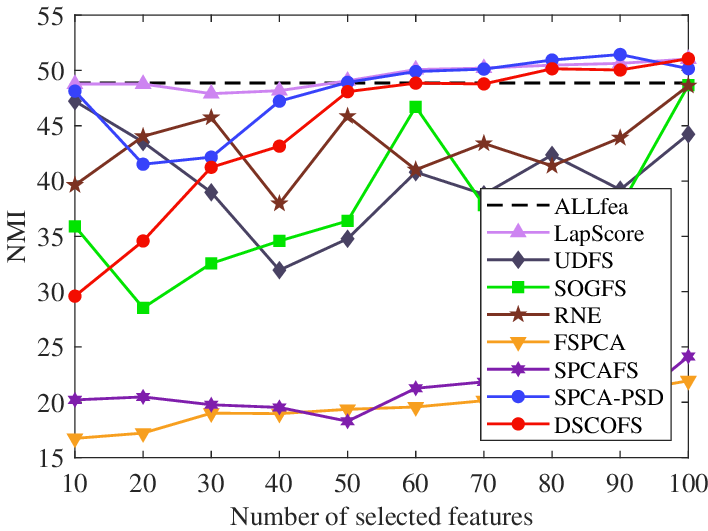}
	}
	\hspace{-4.5mm}
	
	\hspace{-4.5mm}
	\subfigcapskip=-5pt
	\subfigure[UMIST]{
		\label{e}
		\centering
		\includegraphics[width=0.23\textwidth]{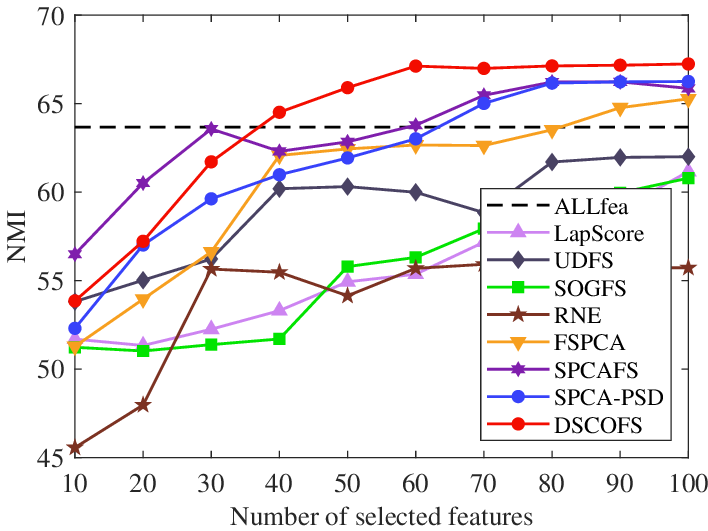}
	}
	\hspace{-1.5mm}  
	\subfigcapskip=-5pt
	\subfigure[warpPIE10P]{
		\label{f}
		\centering
		\includegraphics[width=0.23\textwidth]{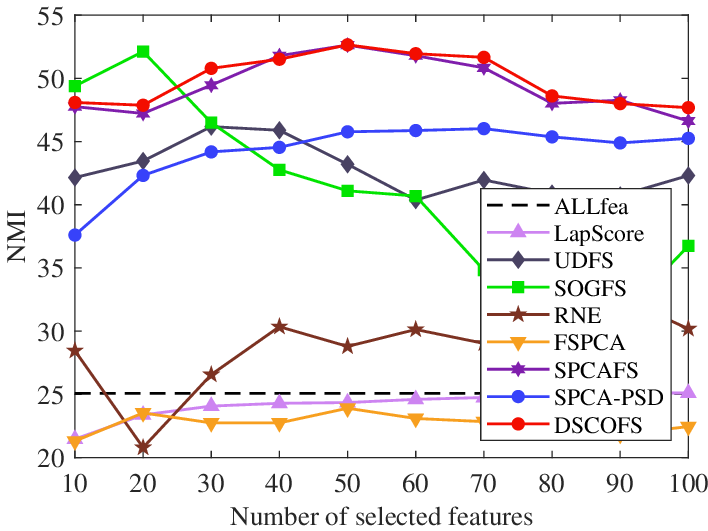}
	}
	\hspace{-1.5mm} 
	\subfigcapskip=-5pt
	\subfigure[Isolet]{
		\label{g}
		\centering
		\includegraphics[width=0.23\textwidth]{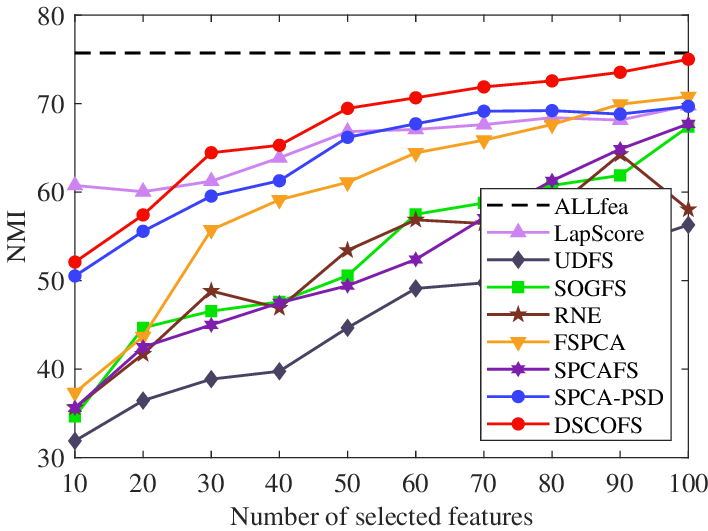}
	}
	\hspace{-1.5mm} 
	\subfigcapskip=-5pt
	\subfigure[MSTAR\_SOC\_CNN]{
		\label{g}
		\centering
		\includegraphics[width=0.23\textwidth]{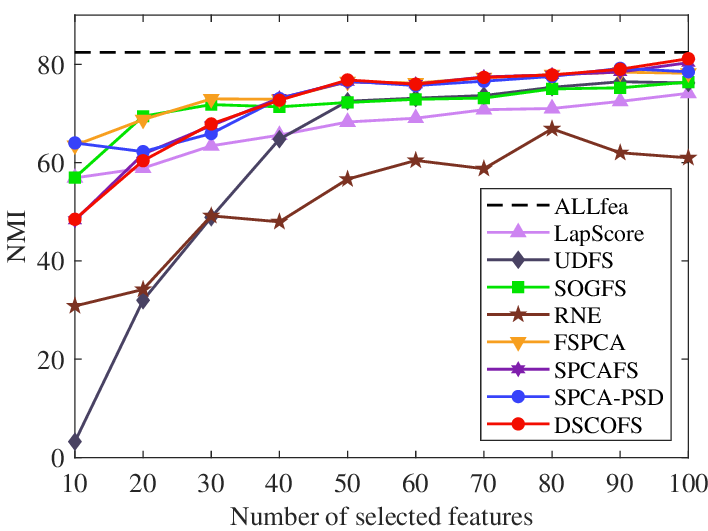}
	}
	\hspace{-4.5mm}
	\vspace{-0.2cm}
	\caption{The NMI (\%) curves of compared methods on eight real-world datasets. All methods select features based on the given parameters and only the mean results of NMI (\%) are plotted.}
	\centering
	\label{plot-nmi}
\end{figure}

\begin{table}[t]
	\caption{The NMI (mean\% $\pm$ std\%) results  of compared methods on eight real-world datasets. The best and second-best results, except ALLfea, are marked in \textcolor[RGB]{255,0,0}{red} and \textcolor[RGB]{0,0,255}{blue}, respectively.}
	\label{nmi}
	\vspace{0.2cm}
	\centering
	\scalebox{0.55}{
		\begin{tabular}{|c|c|c|c|c|c|c|c|c|c|}
			\hline
			Datasets& ALLfea&LapScore&UDFS&SOGFS&RNE&FSPCA&SPCAFS&SPCA-PSD&DSCOFS \\
			\hline\hline
			\multirow{2}*{COIL20}&
			\multirow{2}*{75.37$\pm$1.96}         
			& 69.59$\pm$1.48 & \textcolor[RGB]{0,0,255}{73.54$\pm$1.76} & 68.92$\pm$1.84 & 70.43$\pm$1.92 & 68.50$\pm$1.56 & 69.98$\pm$1.45 & 69.85$\pm$1.41& \textcolor[RGB]{255,0,0}{76.25$\pm$1.71}\\
			&& (100) & \textcolor[RGB]{0,0,255}{(100)} & (100) & (100) & (100) & (100) & (100) & \textcolor[RGB]{255,0,0}{(100)}\\ 
			\multirow{2}*{USPS}&
			\multirow{2}*{61.12$\pm$2.01}          
			& 59.46$\pm$1.80& 54.69$\pm$2.11 & 52.96$\pm$1.54 & 45.36$\pm$1.93 & \textcolor[RGB]{0,0,255}{62.00$\pm$1.87} & 60.98$\pm$2.37 & 60.90$\pm$2.02 & \textcolor[RGB]{255,0,0}{64.06$\pm$2.58}\\
			&& (100) & (100) & (100) & (90) & \textcolor[RGB]{0,0,255}{(60)} & (100) & (100) & \textcolor[RGB]{255,0,0}{(100)}\\ 
			\multirow{2}*{lung\_discrete}&
			\multirow{2}*{62.85$\pm$5.13} 
			& 56.79$\pm$3.99 & 64.84$\pm$5.09 & 59.70$\pm$5.24 & 61.63$\pm$5.83 & 58.26$\pm$6.39 & 69.09$\pm$5.61 & \textcolor[RGB]{0,0,255}{70.93$\pm$5.46} & \textcolor[RGB]{255,0,0}{70.98$\pm$7.00}\\
			&& (100) & (100) & (100) & (70) & (40) & (100) & \textcolor[RGB]{0,0,255}{(80)} & \textcolor[RGB]{255,0,0}{(100)}\\ 
			\multirow{2}*{GLIOMA}&
			\multirow{2}*{48.86$\pm$5.72}
			& 51.03$\pm$2.48 & 47.22$\pm$3.53 & 48.67$\pm$10.98 & 48.62$\pm$6.32 & 21.94$\pm$5.28 & 24.14$\pm$6.97 & \textcolor[RGB]{255,0,0}{51.44$\pm$5.62} & \textcolor[RGB]{0,0,255}{51.06$\pm$6.19}\\
			&& (100) & (10) & (100) & (100) & (100) & (100) & \textcolor[RGB]{255,0,0}{(90)} & \textcolor[RGB]{0,0,255}{(80)}\\  
			\multirow{2}*{UMIST}&
			\multirow{2}*{63.67$\pm$1.85}           
			& 61.16$\pm$1.71 & 62.00$\pm$1.58 & 60.79$\pm$1.54 & 55.92$\pm$1.57 & 65.27$\pm$1.58 & 66.23$\pm$1.60 & \textcolor[RGB]{0,0,255}{66.25$\pm$1.72} & \textcolor[RGB]{255,0,0}{67.24$\pm$1.85}\\
			&& (100) & (100) & (100) & (70) & (100) & (90) & \textcolor[RGB]{0,0,255}{(100)} & \textcolor[RGB]{255,0,0}{(100)}\\  
			\multirow{2}*{warpPIE10P}&
			\multirow{2}*{25.07$\pm$2.88}
			& 25.13$\pm$1.73 & 46.18$\pm$3.30 & 52.12$\pm$3.25 & 32.67$\pm$3.31 & 23.90$\pm$2.01 & \textcolor[RGB]{0,0,255}{52.63$\pm$3.33} & 46.02$\pm$3.70 & \textcolor[RGB]{255,0,0}{52.65$\pm$3.29}\\
			&& (90) & (20) & (20) & (90) & (50) & \textcolor[RGB]{0,0,255}{(50)} & (70) & \textcolor[RGB]{255,0,0}{(50)}\\  
			\multirow{2}*{Isolet}&
			\multirow{2}*{75.72$\pm$1.70}
			& 69.77$\pm$1.20 & 56.29$\pm$1.11 & 67.40$\pm$1.44 & 64.27$\pm$0.95 & \textcolor[RGB]{0,0,255}{70.79$\pm$1.12} & 67.71$\pm$1.33 &  69.69$\pm$0.80 & \textcolor[RGB]{255,0,0}{75.01$\pm$1.35}\\
			&& (100) & (100) & (100) & (90) & \textcolor[RGB]{0,0,255}{(100)} & (100) & (100) & \textcolor[RGB]{255,0,0}{(100)}\\   
			\multirow{2}*{MSTAR\_SOC\_CNN}&
			\multirow{2}*{82.42$\pm$3.31}
			& 74.10$\pm$1.76 & 76.45$\pm$2.47 & 76.39$\pm$1.70 & 66.87$\pm$1.99 & 78.39$\pm$2.17 & \textcolor[RGB]{0,0,255}{80.33$\pm$2.50} & 79.17$\pm$2.77 & \textcolor[RGB]{255,0,0}{81.14$\pm$3.13}\\
			&& (100) & (90) & (100) & (80) & (90) & \textcolor[RGB]{0,0,255}{(100)} & (90) & \textcolor[RGB]{255,0,0}{(100)}\\     
			\hline
			Average & 61.89$\pm$3.07 & 58.38$\pm$2.02 & 60.15$\pm$2.62 & 60.87$\pm$3.44 & 55.72$\pm$2.98 & 56.13$\pm$2.75 & 61.39$\pm$3.15 & \textcolor[RGB]{0,0,255}{64.28$\pm$2.94} & \textcolor[RGB]{255,0,0}{67.30$\pm$3.39}\\
			\hline
	\end{tabular}}
	
\end{table}

\begin{table}[t]
	\caption{The ACC (\%), NMI (\%), and FSR (\%) of ablation study results on eight real-world datasets.}\label{Ablation}
	\vspace{0.2cm}
	\centering
	\scalebox{0.9}{
		\begin{tabular}{|c|c|c|c|c|}
			\hline
			Datasets& $\|X\|_0\leq s $ & ACC & NMI & FSR\\
			\hline\hline
			\multirow{2}*{COIL20}&
			$\times$ & 60.25$\pm$4.52 &75.89$\pm$1.58 &\multirow{2}*{84}\\
			&$\surd$ & 60.51$\pm$4.42 & 76.25$\pm$1.71 &\\ 
			\multirow{2}*{USPS}&
			$\times$ & 67.84$\pm$3.71 & 60.90$\pm$1.95 &\multirow{2}*{68}\\
			&$\surd$ & 69.67$\pm$4.97 & 64.06$\pm$2.58 &\\ 
			\multirow{2}*{lung\_discrete}&
			$\times$ & 71.42$\pm$7.95 & 69.74$\pm$6.11&\multirow{2}*{92}\\
			&$\surd$ & 73.12$\pm$8.48 & 70.98$\pm$7.00&\\ 
			\multirow{2}*{GLIOMA}&
			$\times$ & 58.24$\pm$5.04 & 49.76$\pm$6.12 &\multirow{2}*{85}\\
			&$\surd$ & 60.88$\pm$6.31 &51.06$\pm$6.19 &\\ 
			\multirow{2}*{UMIST}&
			$\times$ & 47.33$\pm$3.05& 67.44$\pm$1.88 &\multirow{2}*{95}\\
			&$\surd$ & 48.10$\pm$3.01 & 67.24$\pm$1.85 &\\ 
			\multirow{2}*{warpPIE10P}&
			$\times$ & 47.91$\pm$4.99& 51.19$\pm$3.79&\multirow{2}*{89}\\
			&$\surd$ & 49.00$\pm$3.88 & 52.65$\pm$3.29 &\\ 
			\multirow{2}*{Isolet}&
			$\times$ & 57.29$\pm$3.44 & 72.82$\pm$1.87 &\multirow{2}*{52}\\
			&$\surd$ & 59.67$\pm$3.46& 75.01$\pm$1.35 &\\ 
			\multirow{2}*{MSTAR\_SOC\_CNN}&
			$\times$ & 82.06$\pm$6.87 & 81.01$\pm$2.41 &\multirow{2}*{99}\\
			&$\surd$ & 82.59$\pm$7.41 & 81.14$\pm$3.13 &\\
			\hline
	\end{tabular}}
\end{table}

\begin{figure}[t]
	\centering
	\hspace{-4.5mm}
	\subfigcapskip=-5pt
	\subfigure[COIL20 ($\ell_{2,0}$)]{
		\centering
		\label{a}
		\includegraphics[width=0.25\textwidth]{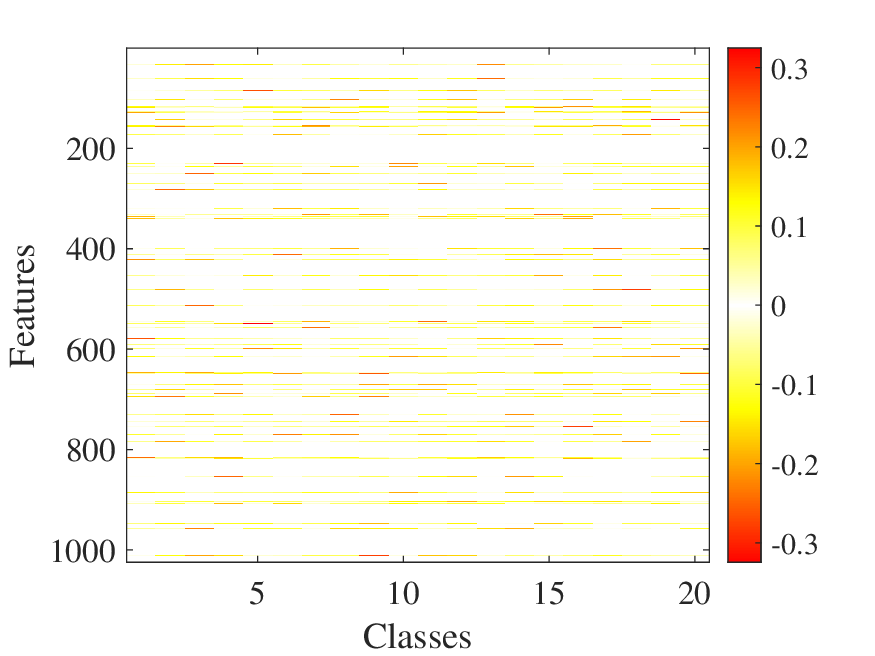}
	}
	\hspace{-5mm}
	\subfigcapskip=-5pt
	\subfigure[USPS ($\ell_{2,0}$)]{
		\centering
		\label{a}
		\includegraphics[width=0.25\textwidth]{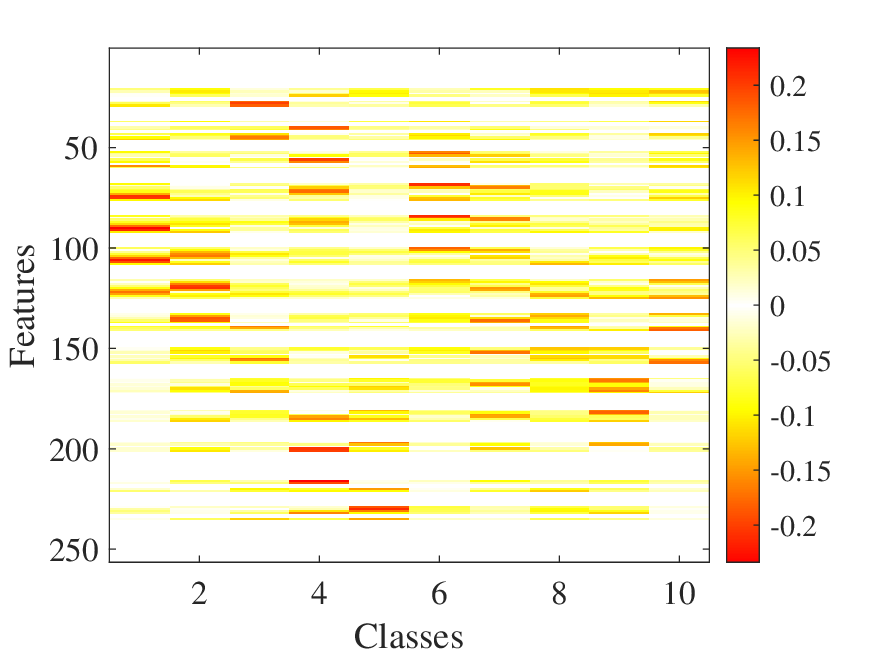}
	}
	\hspace{-5mm} 
	\subfigcapskip=-5pt
	\subfigure[GLIOMA ($\ell_{2,0}$)]{		
		\label{b}
		\centering
		\includegraphics[width=0.25\textwidth]{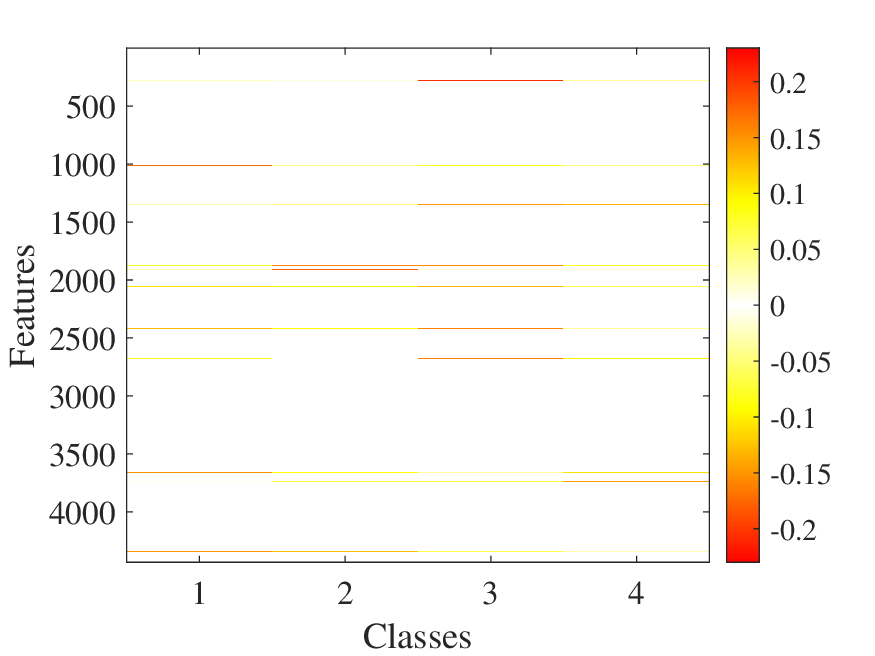}				
	}
	\hspace{-5mm} 
	\subfigcapskip=-5pt
	\subfigure[Isolet ($\ell_{2,0}$)]{		
		\label{b}
		\centering
		\includegraphics[width=0.25\textwidth]{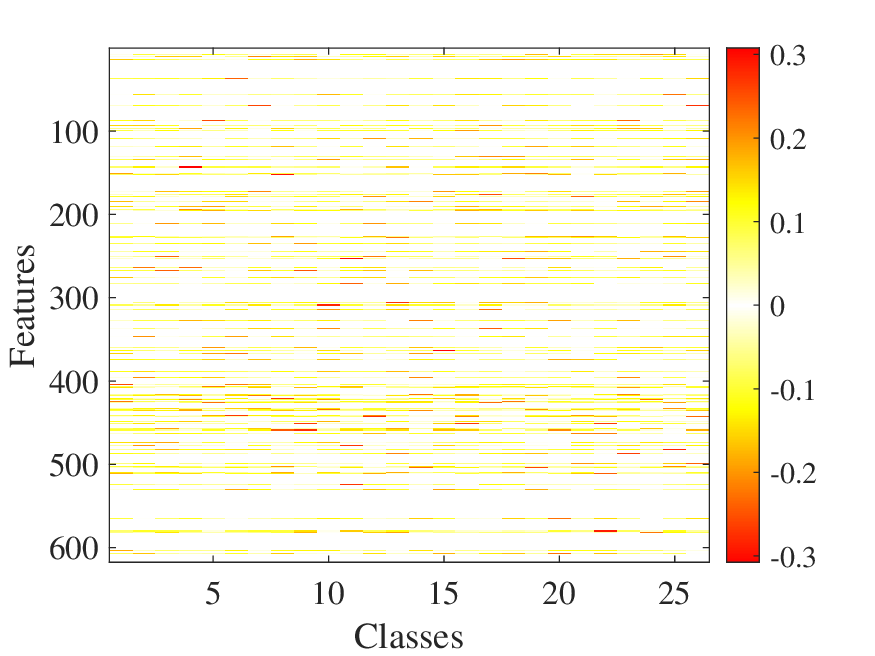}				
	}
	\hspace{-4.5mm}
	
	\hspace{-4.5mm}
	\subfigcapskip=-5pt
	\subfigure[COIL20 ($\ell_{2,0}$+$\ell_{0}$)]{
		\centering
		\label{a}
		\includegraphics[width=0.25\textwidth]{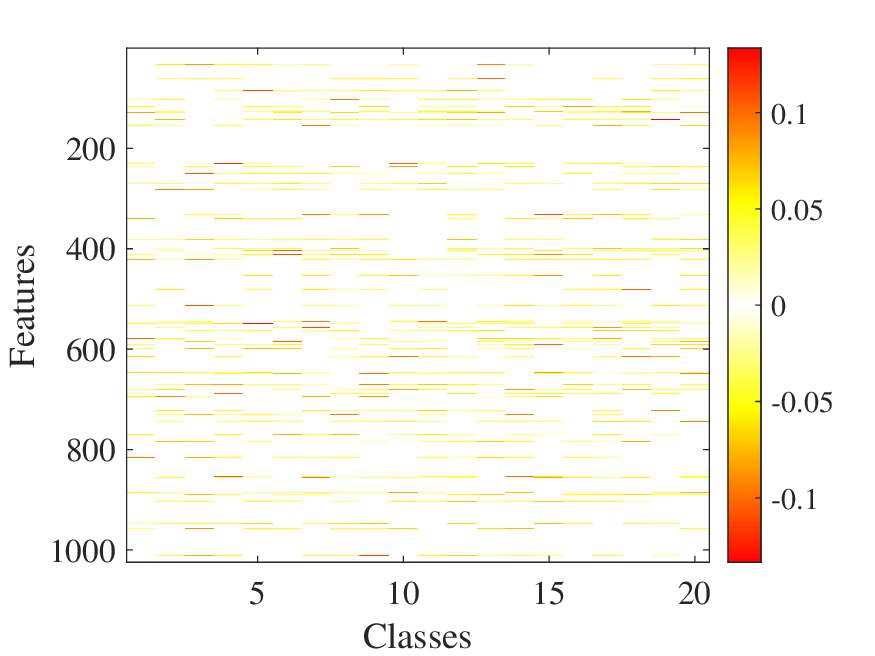}
	}
	\hspace{-5mm}
	\subfigcapskip=-5pt
	\subfigure[USPS ($\ell_{2,0}$+$\ell_{0}$)]{
		\centering
		\label{a}
		\includegraphics[width=0.25\textwidth]{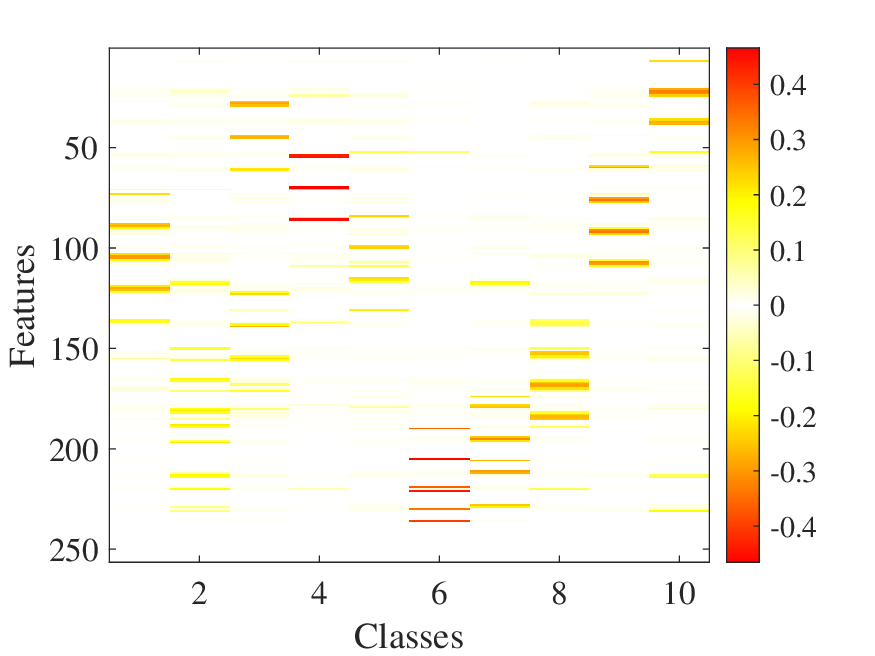}
	}
	\hspace{-5mm} 
	\subfigcapskip=-5pt
	\subfigure[GLIOMA ($\ell_{2,0}$+$\ell_{0}$)]{		
		\label{b}
		\centering
		\includegraphics[width=0.25\textwidth]{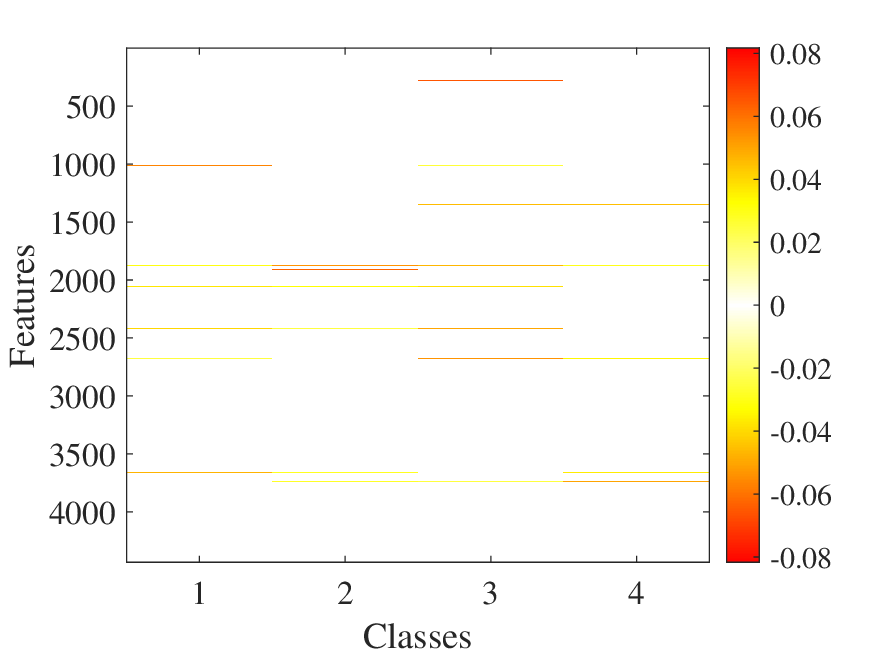}				
	}
	\hspace{-5mm} 
	\subfigcapskip=-5pt
	\subfigure[\small{Isolet ($\ell_{2,0}$+$\ell_{0}$)}]{		
		\label{b}
		\centering
		\includegraphics[width=0.25\textwidth]{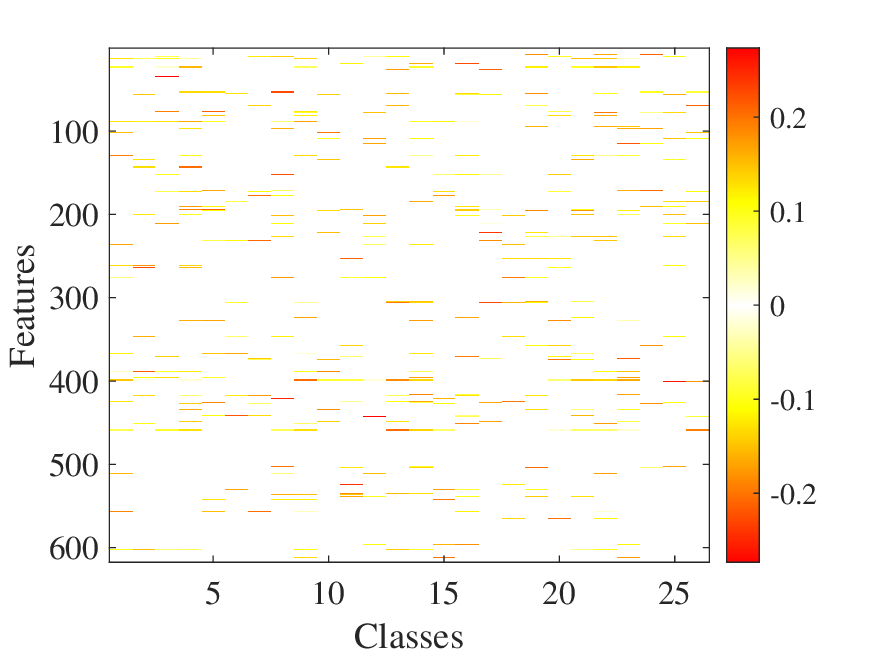}				
	}
	\hspace{-4.5mm}
	\vspace{-0.2cm}
	\caption{The sparse visualization of the projection matrix $X$ on four real-world datasets, where white means the value of the element is 0. (a)-(d) are the $X$ learned from $\ell_{2,0}$-norm; (e)-(h) are from double sparsity.}
	\centering
	\label{plot-heat1}
\end{figure}

\begin{table}[t]
	\caption{The Friedman test results in terms of ACC.}\label{FT}
	\vspace{0.2cm}
	\centering
	\begin{tabular}{|c|c|c|c|}
		\hline
		Methods& Ranking&  $\varrho$-value&Hypothesis\\
		\hline \hline
		LapScore&6&\multirow{8}*{0.0001}&\multirow{8}*{Reject}\\
		UDFS&4.75&&\\
		SOGFS&5.75&&\\
		RNE&6.125&&\\
		FSPCA&5.5&&\\
		SPCAFS&3.75&&\\
		SPCA-PSD&3.125&&\\
		DSCOFS&1&&\\
		\hline
	\end{tabular}
\end{table}

In addition, Figure \ref{plot-nmi} shows the corresponding mean curves of NMI. Obviously, our proposed DSCOFS outperforms the baseline of all datasets except Isolet and MSTAR\_S\\OC\_CNN. Compared with other comparison methods, our proposed DSCOFS achieves the best results on all datasets except GLIOMA.  Table \ref{nmi} lists the corresponding mean, standard deviation, and number of selected features. Obviously, the average NMI of all datasets is improved by 3.02\%. Similar to the ACC results, the NMI results show significant improvements on COIL20 and Isolet. In particular, our proposed DSCOFS is at least 2.71\% gains on COIL20 and 4.22\% on Isolet. It is particularly noteworthy that FSPCA, SPCAFS, SPCA-PSD, and our proposed DSCOFS obtain the two best clustering results on all datasets except COIL20. This shows that SPCA is promising for feature selection, but the clustering results obtained with different sparsity show different performances on different datasets.

In summary, the ACC and NMI results demonstrate that our proposed DSCOFS is effective and performs well on all datasets. 
There comes a conclusion: double sparsity allows us to deal with more complicated structures than single sparsity, which also suggests the potential of our proposed DSCOFS.

\subsection{Analysis of Double Sparsity}\label{exp-3}

Although the above feature selection experiments verify the effectiveness of our proposed DSCOFS, they cannot intuitively illustrate the role of element-wise sparsity. Accordingly, this subsection conducts the ablation study and statistical tests to analyze the effect of element-wise sparsity in feature selection.

\subsubsection{Ablation Study}

For consistency, the same parameter setting rules, feature selection rules, and initial solutions are used here. Table \ref{Ablation} records the best ACC and the corresponding NMI within 100 features, where $\times$ means that only $\ell_{2,0}$-norm structure sparsity is considered, and $\surd$ means that double sparsity is considered.

In order to reflect the differences of different methods, we introduce a new evaluation metric called feature similarity rate (FSR). 
To the best of our knowledge, this metric has never been considered in previous work on feature selection.
Let $\mathbb{T}_\textrm{DSCOFS}$, $\mathbb{T}_{2,0}$ denote the sets of selected features, where the features are the rankings obtained by executing double sparsity in our proposed DSCOFS and only $\ell_{2,0}$-norm, respectively. FSR is defined as
\begin{equation}\label{fsr}
\begin{aligned}
\textrm{FSR}=\frac{1}{n} \textrm{card}(\mathbb{T}_\textrm{DSCOFS}\cap \mathbb{T}_{2,0}),
\end{aligned}
\end{equation}
where $n$ is the number of features selected when calculating FSR. In this study, $n$ is fixed to 100.  Naturally, FSR indicates the percentage of features in two sets of features election results that are overlapped. The smaller the FSR, the greater the difference.

As can be seen from Table \ref{Ablation}, by introducing element-wise sparsity, the ACC and NMI of our proposed DSCOFS have been improved to a certain extent. The FSR index varies greatly on different datasets. Especially on USPS and Isolet, their FSR values are 68\% and 52\%, respectively, indicating that adding element-wise sparsity selects different features than those selected by only structural sparsity. Furthermore, Figure \ref{plot-heat1} visualizes the projection matrix on four real-world datasets including COIL20, USPS, GLIOMA, and Isolet. It can be found that the sparse transformation matrices (see (e)-(h) in Figure \ref{plot-heat1}) obtained using double sparsity retain different features from the sparse transformation matrices (see (a)-(d) in Figure \ref{plot-heat1}) obtained using only $\ell_{2,0}$-norm. The difference is more evident for USPS and Isolet, which explains why our proposed DSCOFS has a significant improvement on these two datasets.

\subsubsection{Statistical Tests}

The Friedman test is a ranking-based statistical method used to compare whether there is a significant difference in the average performance of multiple methods. Initially, the ACC of each method on each dataset is ranked from best to worst, assigning rank values from 1 to 8. When more than one method has the same accuracy on the same dataset, the average ranking values are taken. In this experiment, the null hypothesis $\mathcal{H}_0$ of the Friedman test is that there is no significant difference in the performance of all comparison methods. Then, with a significance level set at $\alpha = 0.05$, it can be seen from Table \ref{FT} that $\varrho=0.0001$, which means that the result rejects the null hypothesis $\mathcal{H}_0$. This convinces us to believe that there exist significant differences among all comparison methods.

\begin{figure}[t]
	\centering
	\includegraphics[width=0.80\textwidth]{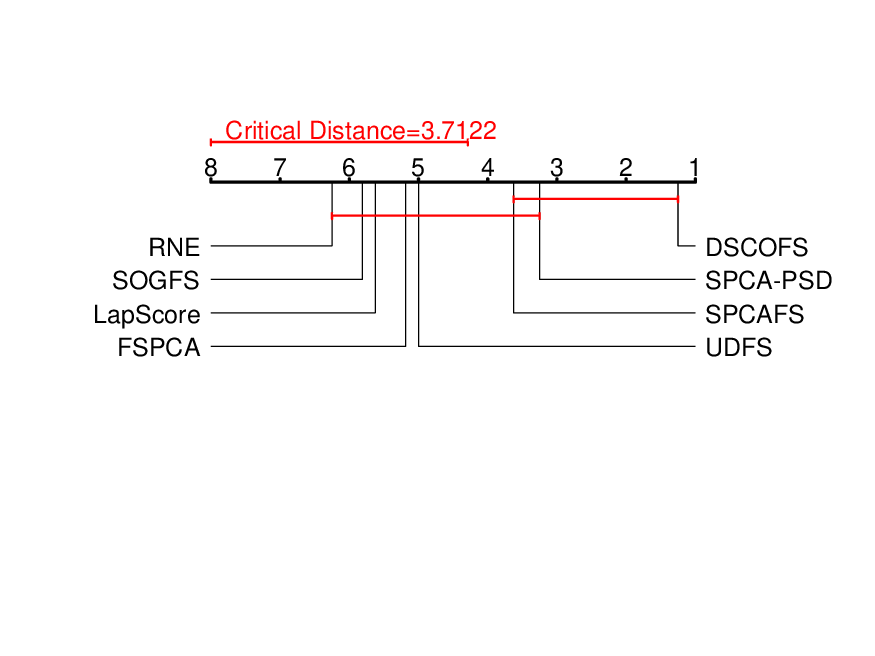}
	\vspace{-3.5cm}
	\caption{The post-hoc Nemenyi test results in terms of ACC.}
	\label{ft}
\end{figure}

However, the Friedman test cannot determine whether there is a difference between the two methods. This is what the post-hoc Nemenyi test is intended to solve, which can measure the difference through the critical difference (CD) value. Figure \ref{ft} shows the corresponding results of the post-hoc Nemenyi test. Our proposed DSCOFS is significantly different from LapScore, UDFS, SOGFS, RNE, and FSPCA, but not significantly different from SPCA-PSD and SPCAFS. Even so, DSCOFS has better performance than them. In addition, it is worth noting that there is a significant difference between our proposed DSCOFS and FSPCA with $\ell_{2,0}$-norm constraint, which also indicates the advantages of double sparsity. 

\subsection{Discussion}\label{exp-4}

After experimentally validating and analyzing the effectiveness of our proposed DSCOFS, this subsection further discusses the 
parameter sensitivity, model stability, and convergence on four datasets. In addition, the comparison with deep learning-based UDS methods is analyzed.

\begin{figure}[t]
	\centering
	\hspace{-4.5mm}
	\subfigcapskip=-5pt
	\subfigure[COIL20 (ACC)]{
		\centering
		\label{a}
		\includegraphics[width=0.26\textwidth]{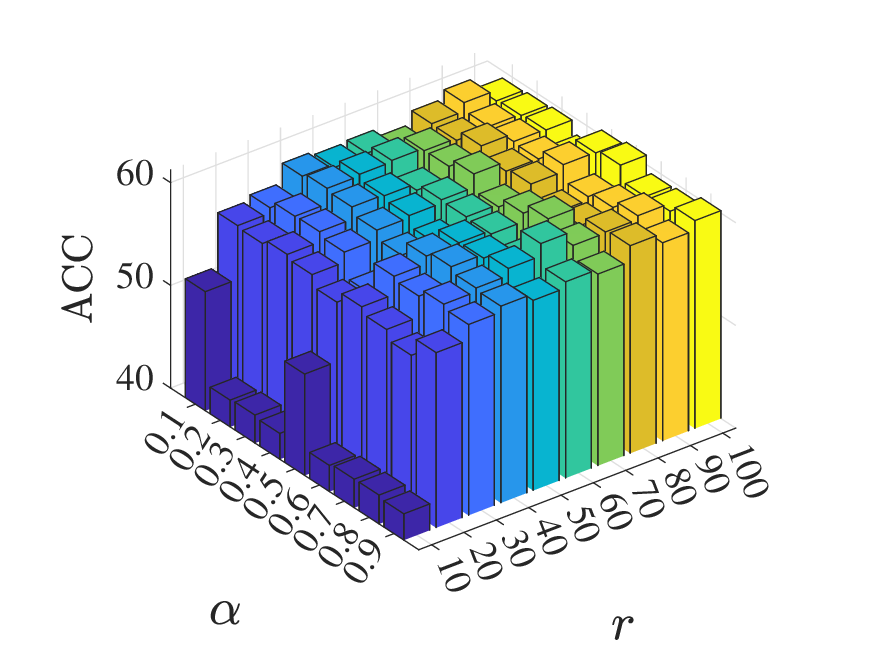}
	}
	\hspace{-8mm} 
	\subfigcapskip=-5pt
	\subfigure[USPS (ACC)]{		
		\label{b}
		\centering
		\includegraphics[width=0.26\textwidth]{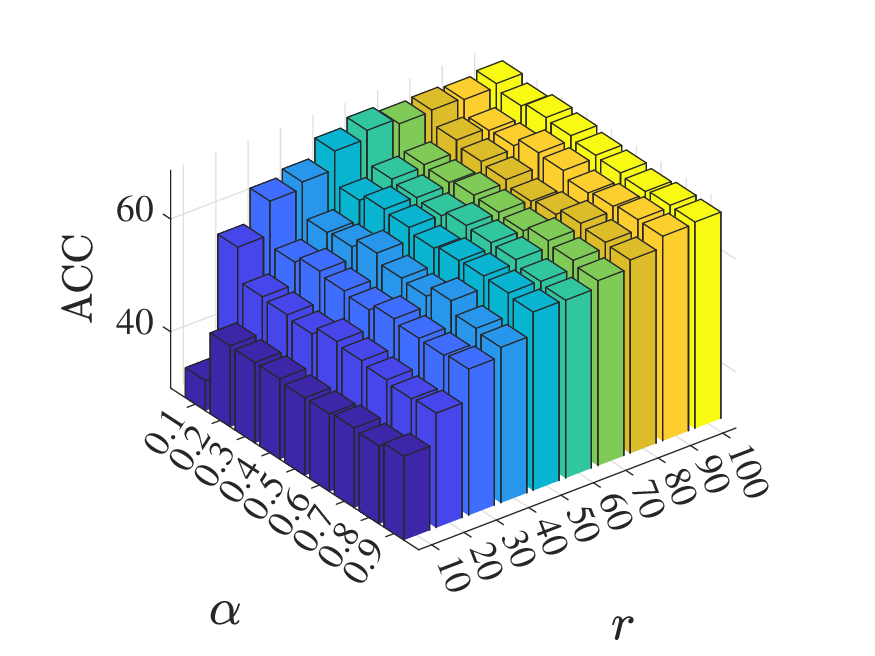}				
	}
	\hspace{-8mm}
	\subfigcapskip=-5pt
	\subfigure[lung\_discrete (ACC)]{
		\label{c}
		\centering
		\includegraphics[width=0.26\textwidth]{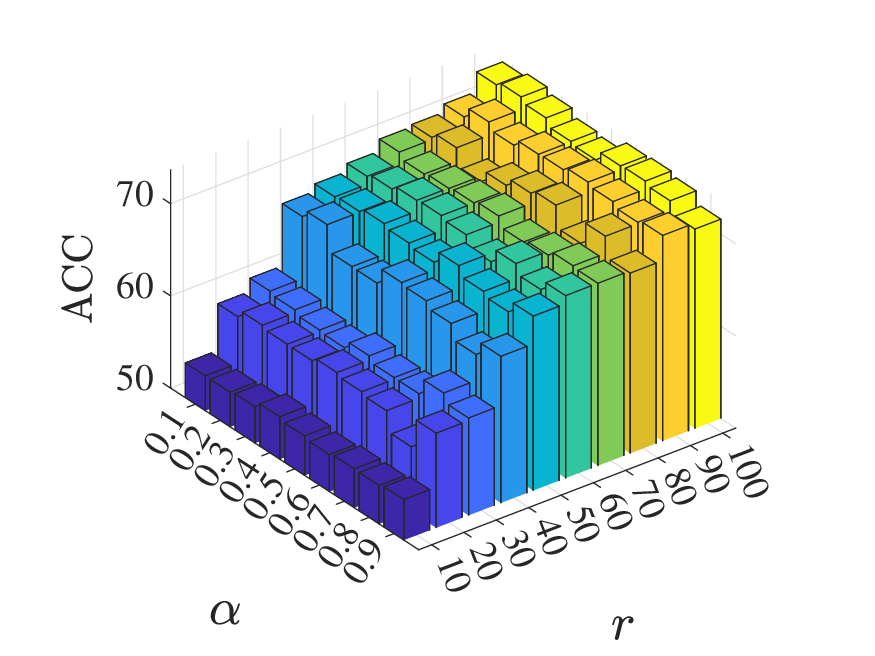}
	}
	\hspace{-8mm} 
	\subfigcapskip=-5pt
	\subfigure[Isolet (ACC)]{
		\label{d}
		\centering
		\includegraphics[width=0.26\textwidth]{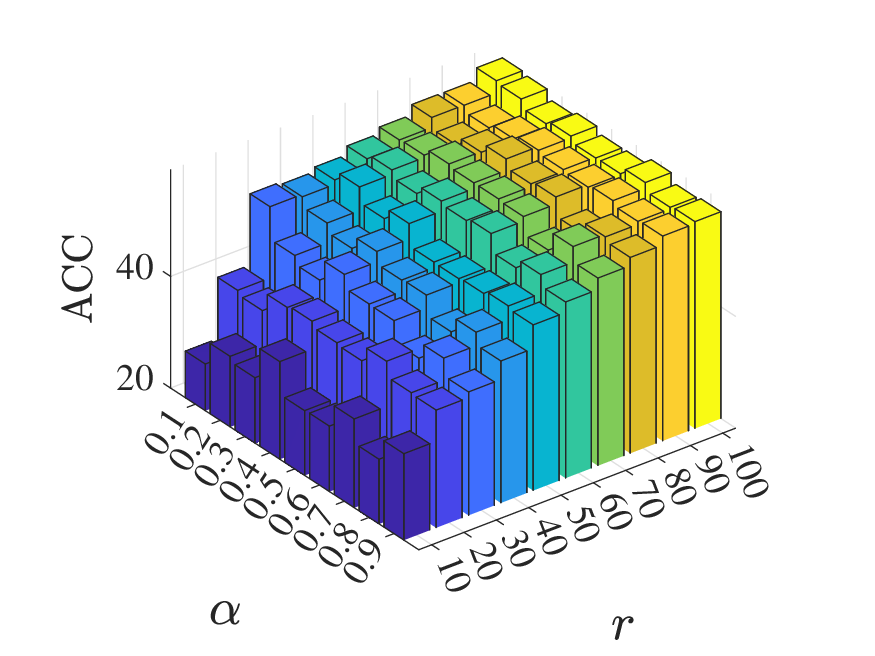}
	}
	\hspace{-4.5mm} 
	
	\vspace{-2.5mm}   
	\hspace{-4.5mm}
	\subfigcapskip=-5pt
	\subfigure[COIL20 (NMI)]{
		\centering
		\label{a}
		\includegraphics[width=0.26\textwidth]{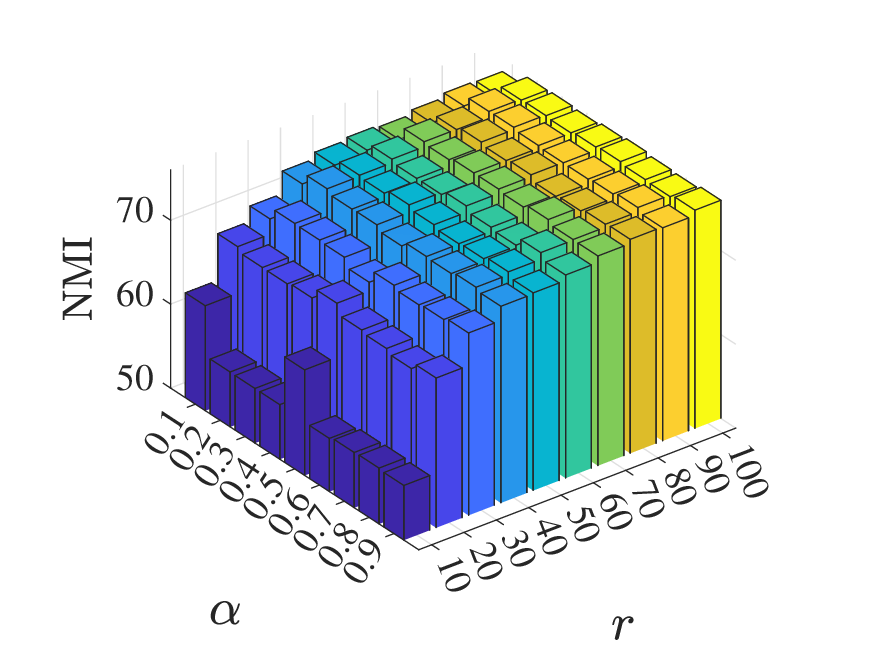}
	}
	\hspace{-8mm} 
	\subfigcapskip=-5pt
	\subfigure[USPS (NMI)]{		
		\label{b}
		\centering
		\includegraphics[width=0.26\textwidth]{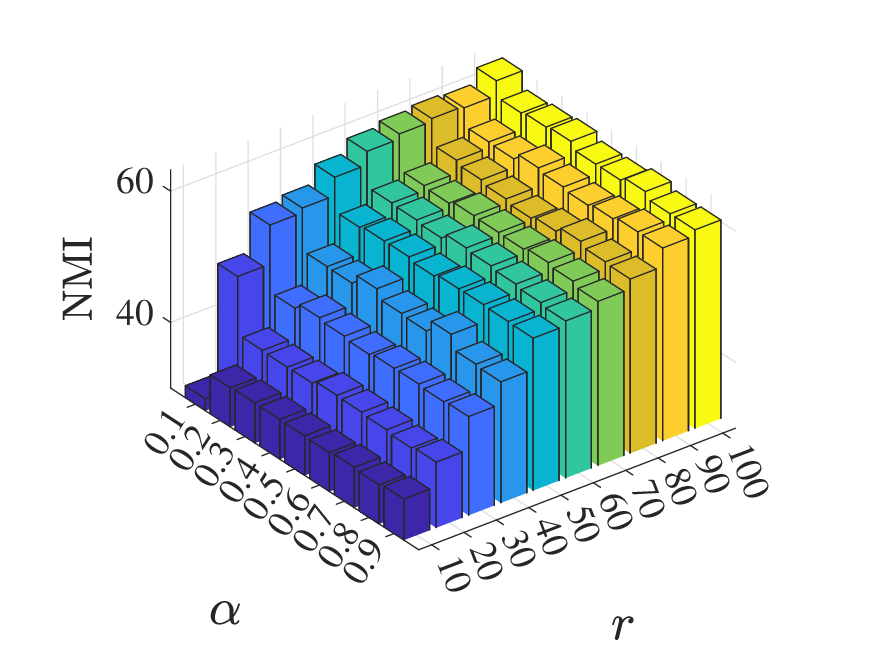}				
	}
	\hspace{-8mm}
	\subfigcapskip=-5pt
	\subfigure[lung\_discrete (NMI)]{
		\label{c}
		\centering
		\includegraphics[width=0.26\textwidth]{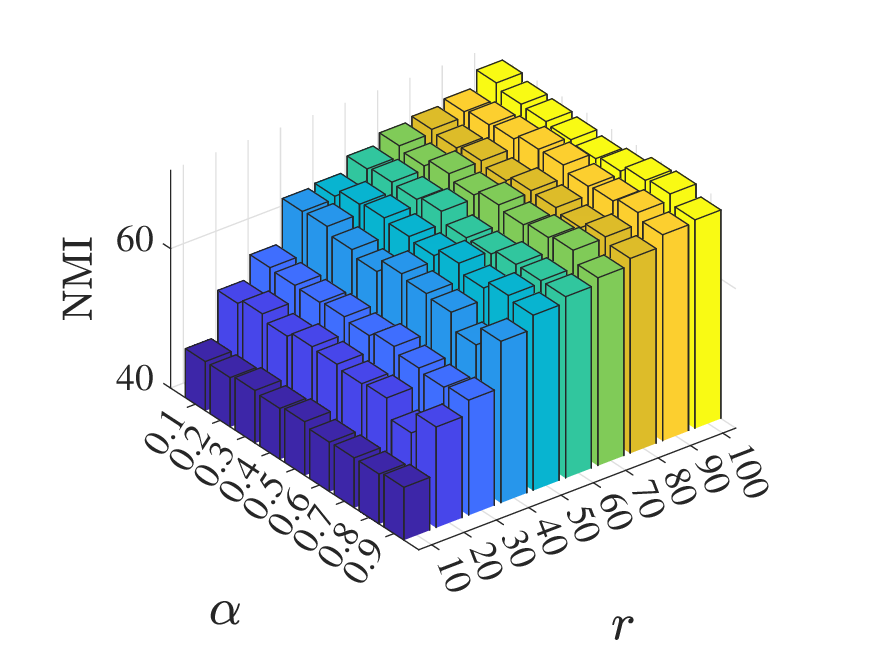}
	}
	\hspace{-8mm} 
	\subfigcapskip=-5pt
	\subfigure[Isolet (NMI)]{
		\label{d}
		\centering
		\includegraphics[width=0.26\textwidth]{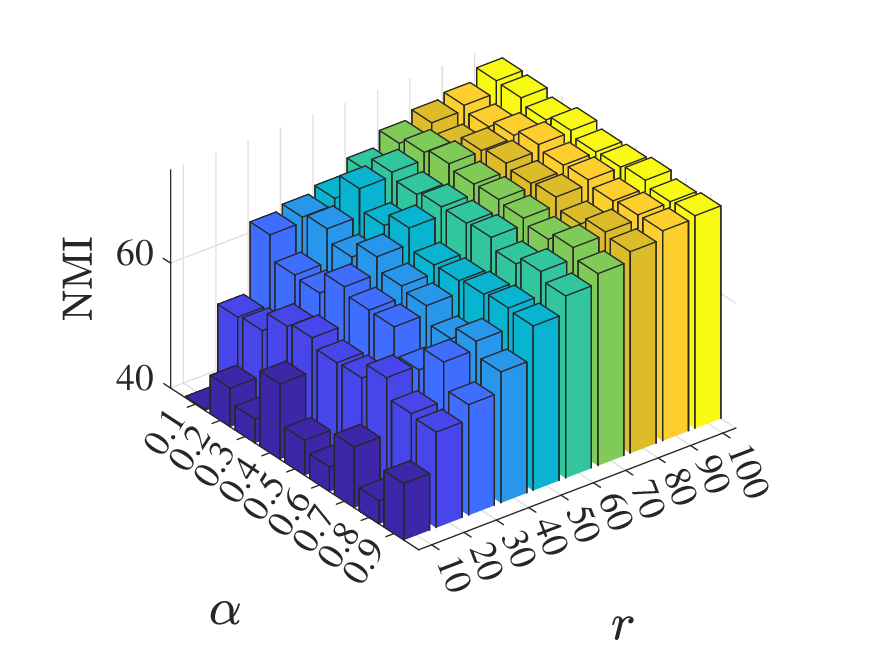}
	}
	\hspace{-4.5mm} 
	\vspace{-0.2cm}   
	\caption{The ACC (\%) and NMI (\%) with different values of our proposed DSCOFS on four real-world datasets. (a)-(d) are the ACC (\%) results; (e)-(h) are the NMI (\%) results.}
	\centering
	\label{plot-psa}
\end{figure}

\begin{figure}[!htp]
	\centering
	\hspace{-4.5mm}
	\subfigcapskip=-5pt
	\subfigure[COIL20 (ACC)]{
		\centering
		\label{a}
		\includegraphics[width=0.25\textwidth]{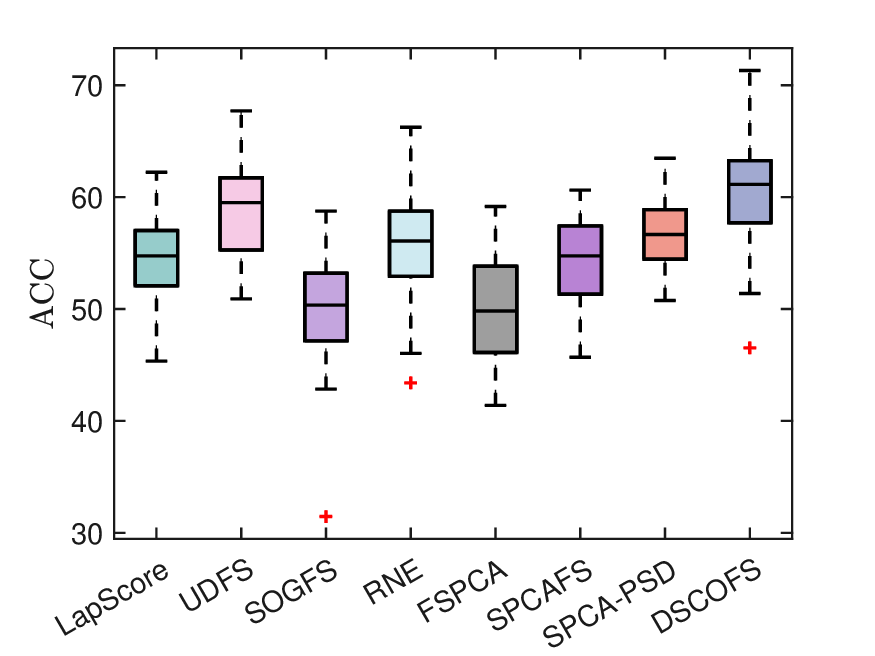}
	}
	\hspace{-4.5mm} 
	\subfigcapskip=-5pt
	\subfigure[USPS (ACC)]{		
		\label{b}
		\centering
		\includegraphics[width=0.25\textwidth]{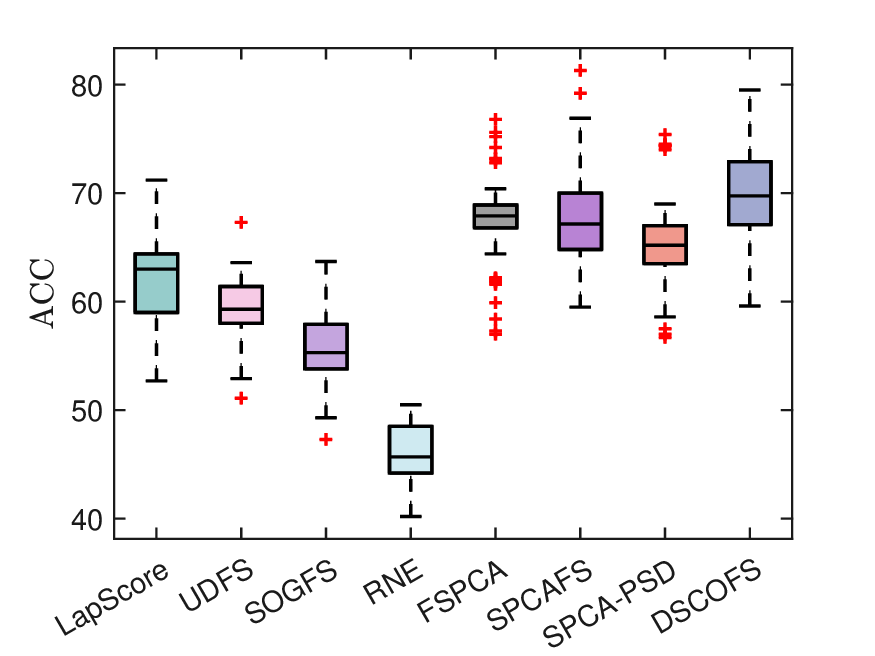}				
	}
	\hspace{-4.5mm} 
	\subfigcapskip=-5pt
	\subfigure[lung\_discrete (ACC)]{
		\centering
		\label{a}
		\includegraphics[width=0.25\textwidth]{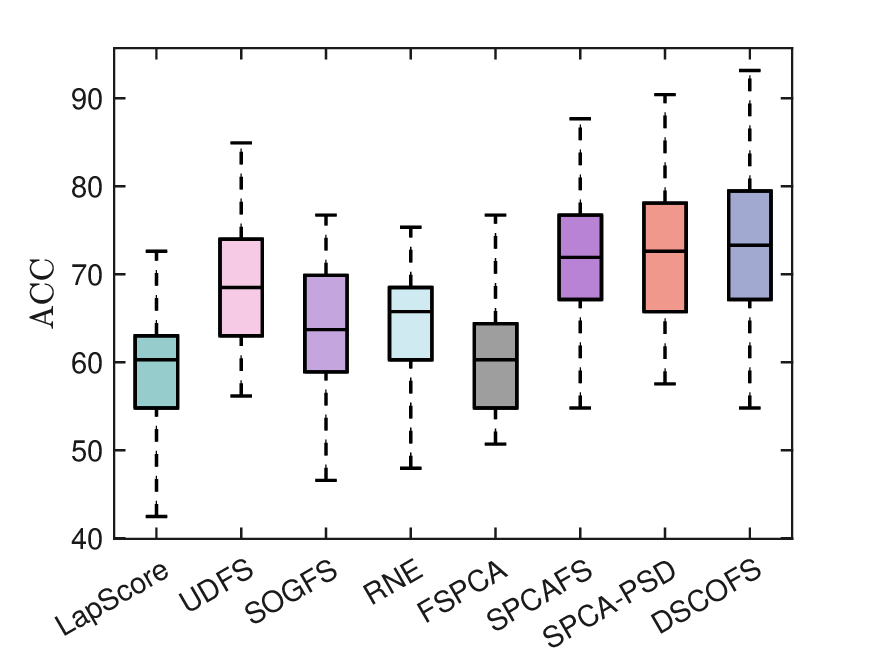}
	}
	\hspace{-4.5mm} 
	\subfigcapskip=-5pt
	\subfigure[Isolet (ACC)]{		
		\label{b}
		\centering
		\includegraphics[width=0.25\textwidth]{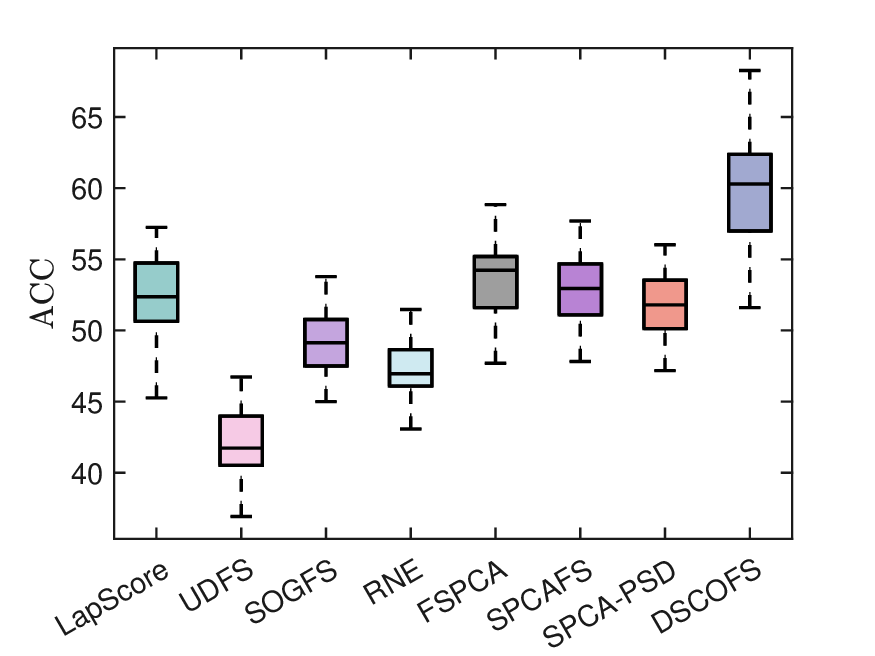}				
	}
	\hspace{-4.5mm}
	
	\vspace{-2.5mm}
	\hspace{-4.5mm}
	\subfigcapskip=-5pt
	\subfigure[COIL20 (NMI)]{
		\centering
		\label{a}
		\includegraphics[width=0.25\textwidth]{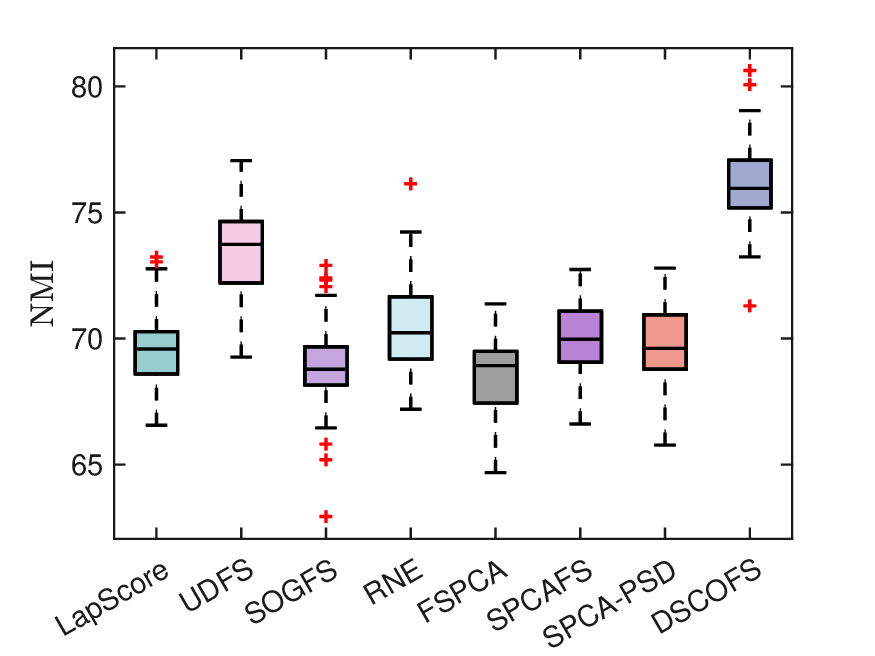}
	}
	\hspace{-4.5mm} 
	\subfigcapskip=-5pt
	\subfigure[USPS (NMI)]{		
		\label{b}
		\centering
		\includegraphics[width=0.25\textwidth]{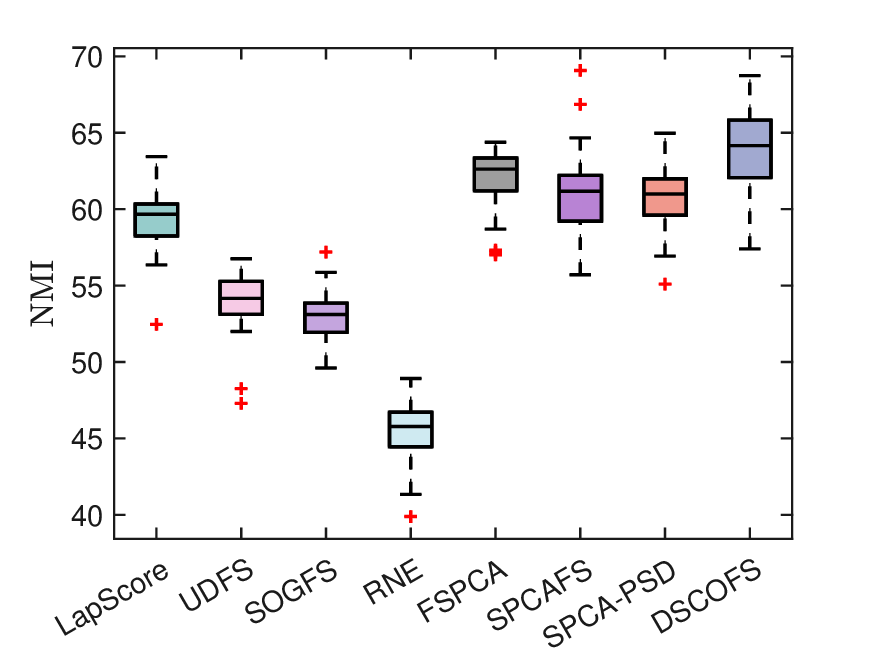}				
	}
	\hspace{-4.5mm} 
	\subfigcapskip=-5pt
	\subfigure[lung\_discrete (NMI)]{
		\centering
		\label{a}
		\includegraphics[width=0.25\textwidth]{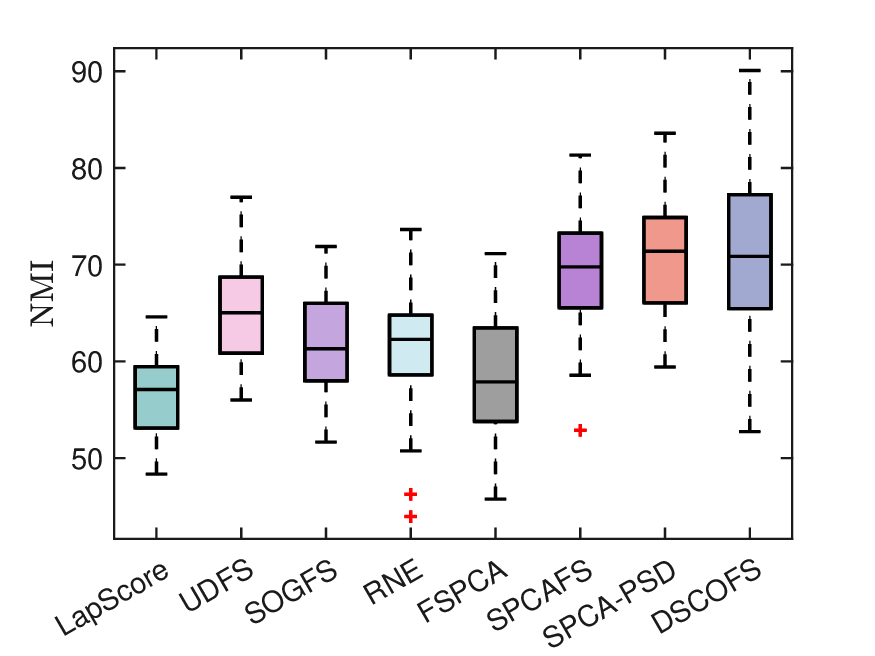}
	}
	\hspace{-4.5mm} 
	\subfigcapskip=-5pt
	\subfigure[Isolet (NMI)]{		
		\label{b}
		\centering
		\includegraphics[width=0.25\textwidth]{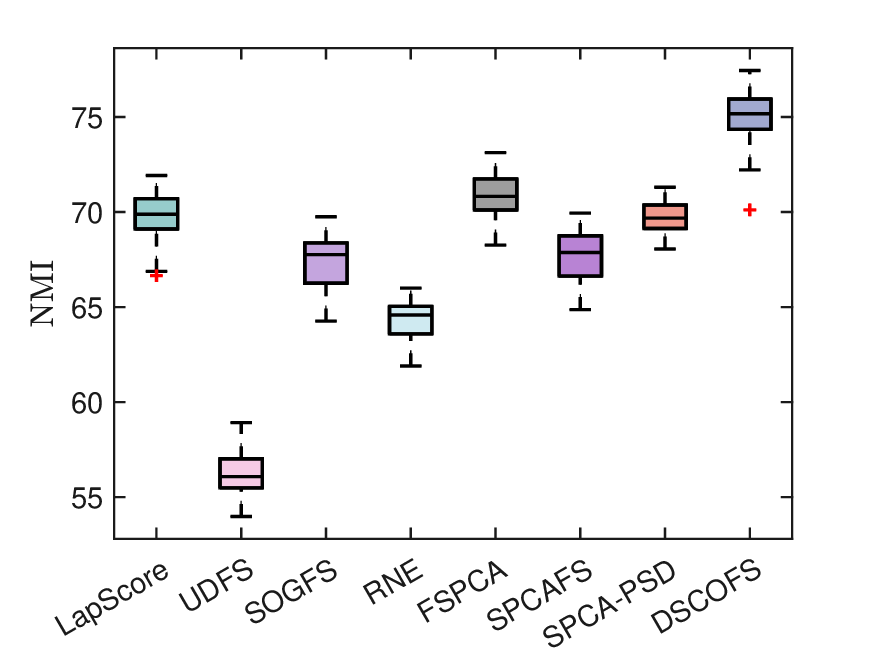}				
	}
	\hspace{-4.5mm}
	\vspace{-0.2cm}   
	\caption{The model stability of ACC(\%) and NMI(\%) of all compared methods on four real-world datasets. (a)-(d) are the ACC (\%) results; (e)-(h) are the NMI (\%) results.}
	\centering
	\label{plot-box}
\end{figure}

\begin{figure}[t]
	\centering
	\hspace{-4.5mm}
	\subfigcapskip=-5pt
	\subfigure[COIL20]{
		\centering
		\label{a}
		\includegraphics[width=0.32\textwidth]{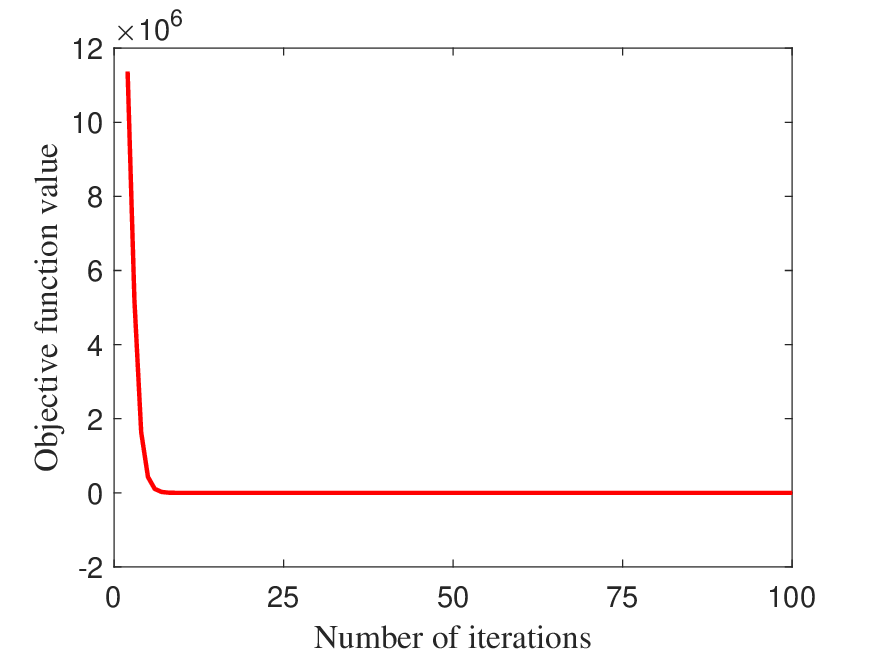}
	}
	\hspace{-2mm} 
	\subfigcapskip=-5pt
	\subfigure[USPS]{		
		\label{b}
		\centering
		\includegraphics[width=0.32\textwidth]{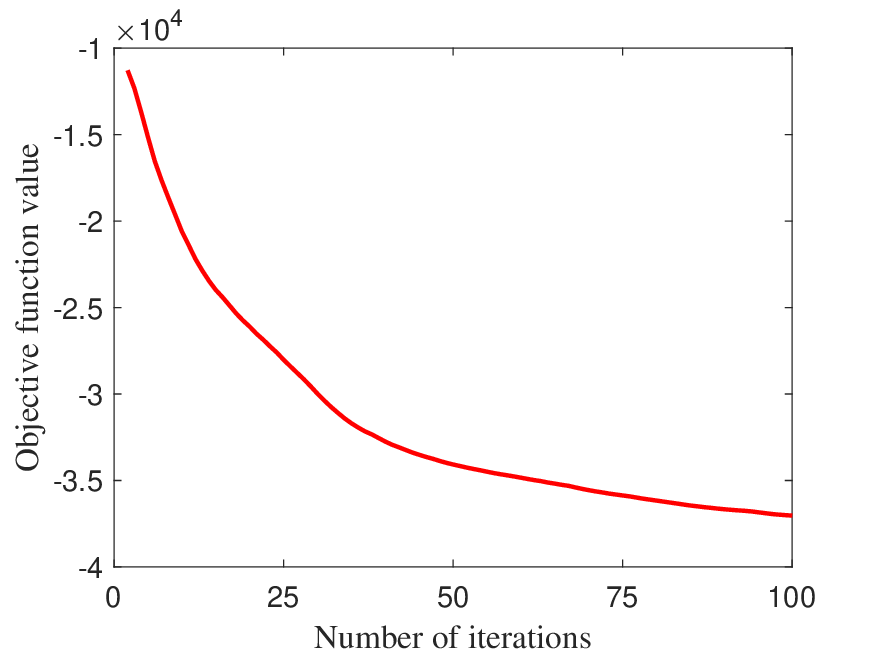}				
	}
	\hspace{-4.5mm}
	
	\vspace{-2.5mm}  
	\hspace{-4.5mm}
	\subfigcapskip=-5pt
	\subfigure[lung\_discrete]{
		\centering
		\label{a}
		\includegraphics[width=0.32\textwidth]{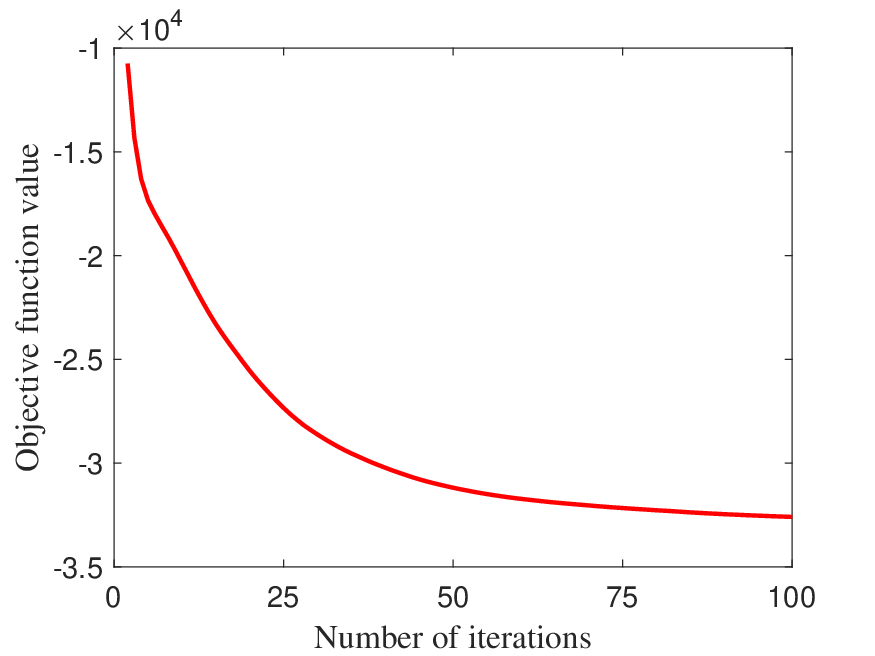}
	}
	\hspace{-2mm} 
	\subfigcapskip=-5pt
	\subfigure[Isolet]{		
		\label{b}
		\centering
		\includegraphics[width=0.32\textwidth]{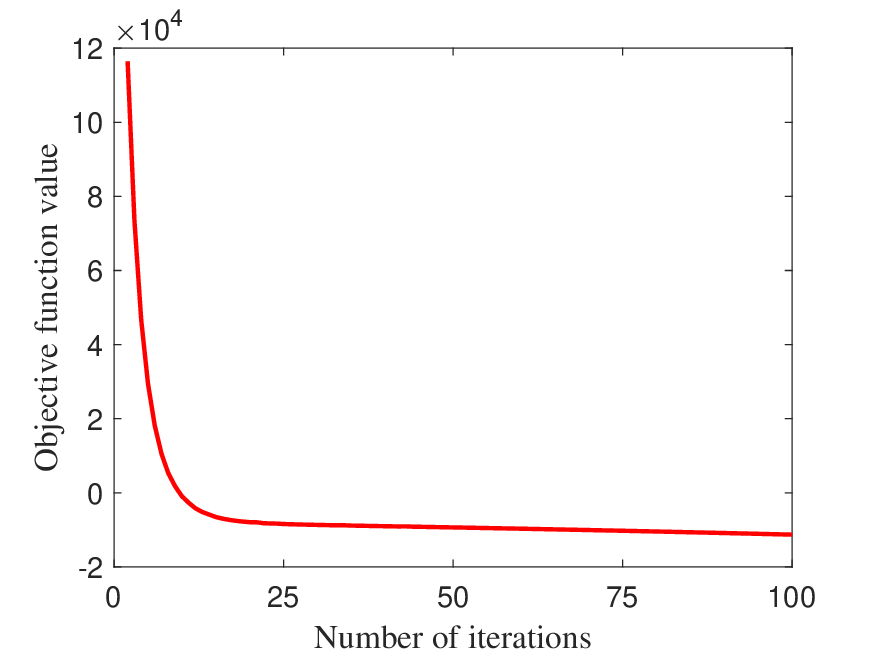}				
	}
	\hspace{-4.5mm}
	\vspace{-0.2cm}   
	\caption{The convergence curves of our proposed DSCOFS on four real-world datasets.}
	\centering
	\label{plot-con}
\end{figure}

\subsubsection{Parameter Sensitivity} \label{para}

Figure \ref{plot-psa} shows the best clustering results under different $r$ and $\alpha$, where $\alpha$ represents the sparsity percentage. It can be seen that the alterations in ACC and NMI are significant under different element-wise sparsity percentages. For USPS, it has a significant performance improvement when $\alpha=0.1$. This also implies that the element-wise  sparsity enhances the model's capacity for feature selection.

\subsubsection{Model Stability}

Figure \ref{plot-box} shows the box plots of the mean and standard deviation of 50 clustering results. It can be observed that the clustering results fluctuate significantly, however, the overall results of our proposed DSCOFS are better than those of other comparison methods. Especially on Isolet, the lowest value  exceeds or is close to the median value of other methods. Although the maximum and minimum values of our proposed DSCOFS on lung\_discrete vary greatly, the average result is still slightly better than other methods. All these illustrate the stability of our proposed method.

\subsubsection{Convergence Analysis}

Figure \ref{plot-con} shows the convergence curves associated with the objective function of Algorithm \ref{am1} when performing feature selection. The results demonstrate that our proposed DSCOFS can continue to decrease and reach a stable state within 100 iterations in most cases. This is consistent with the conclusion in Section \ref{convergence}.

\subsubsection{Comparsion with Deep Learning}

Table \ref{comparison} presents a comprehensive comparison between our proposed DSCOFS with a deep learning-based method, i.e., TSFS+TSNE \cite{mirzaei2020deep}. 
It shows that DSCOFS is comparable to TSFS+TSNE and even outperforms it on the USPS and lung\_discrete datasets, with 7.71\% and 8.76\% higher accuracy, respectively. 
This may be because deep learning-based UFS methods like TSFS+TSNE require large sample datasets for training, 
but our proposed DSCOFS performs well on smaller datasets due to its enhanced interpretability and generalization capabilities. 

\begin{table}[t] 
	\caption{The ACC (\%) and NMI (\%) results of DSCOFS and TSFS+TSNE on four real-world datasets.}\label{comparison}
	\vspace{0.2cm}   
	\centering
	\begin{tabular}{|c|c|c|c|}
		\hline
		Datasets&  Methods & ACC & NMI \\
		\hline\hline
		\multirow{2}*{COIL20}&
		TSFS+TSNE    &{60.80$\pm$3.83} & {71.59$\pm$1.46} \\
		&DSCOFS & {60.51$\pm$4.42} & {76.25$\pm$1.71} \\ 
		\multirow{2}*{USPS}&
		TSFS+TSNE & {61.96$\pm$3.96} & {56.20$\pm$1.20} \\
		&DSCOFS & {69.67$\pm$4.97} & {64.06$\pm$2.58} \\ 
		\multirow{2}*{lung\_discrete}&
		TSFS+TSNE & {64.36$\pm$7.24} & {61.61$\pm$5.70} \\
		&DSCOFS & {73.12$\pm$8.48} & {70.98$\pm$7.00} \\ 
		\multirow{2}*{Isolet}&
		TSFS+TSNE & {60.40$\pm$4.34} & {76.13$\pm$1.54} \\
		&DSCOFS & {59.67$\pm$3.46} & {75.01$\pm$1.35} \\ 
		\hline
	\end{tabular}
\end{table}

\section{Conclusion} \label{conclusion}

In this paper, we have proposed a novel double sparsity constrained PCA model that can select informative features from unlabeled datasets by integrating $\ell_{2,0}$-norm and $\ell_0$-norm constraints. By leveraging the double sparsity, we can promote sparsity in the transformation matrix for filtering out redundant and irrelevant features, thereby achieving more accurate and effective feature selection. In algorithms, we have developed a convergent optimization scheme that incorporates an exact penalty function method into the proximal alternating minimization framework. Extensive numerical experiments have validated the superiority of our proposed method. In particular, we have introduced a new metric, which intuitively analyzes that element-wise sparsity enhances the feature selection ability and addresses the shortcomings of single structural sparsity. It is worth noting that although this work mainly focuses on feature selection, double sparsity constrained optimization can be easily extended to other tasks in image and video processing.

Of course, there are still some limitations to be investigated in the future. On the one hand, it is necessary to develop second-order or distributed optimization algorithms when dealing with the computational challenges brought by high-dimensional data. On the other hand, it is possible to consider using deep unfolding networks to develop model-data-driven methods to learn structures and avoid parameter selection.

\section*{Acknowledgements}

This work was supported in part by the National Natural Science Foundation of China under Grant 12371306, the Innovation Program of Shanghai Municipal Education Commission under Grant 2023ZKZD47, and the 111 Project under Grant B16009.

\section*{Appendix}
\noindent{\textbf{Proof of Theorem \ref{Theor1} (a).}}
Let $X^{k+1}$, $Y^{k+1}$ and $Z^{k+1}$ be the solutions of \eqref{sub-x}, \eqref{sub-y} and \eqref{sub-z}, respectively.  For any points $X^{k}\in\mathcal{M}$, $Y^{k}\in\mathcal{S}$ and $Z^{k}\in\mathcal{R}$, it holds 
\begin{equation}\label{fxyz-x}
\begin{aligned}
f(X^{k+1},Y^k,  Z^k)&\leq f(X^k, Y^k,  Z^k)-\tau_1\|X^{k+1}-X^k\|_\textrm{F}^2,\\
f(X^{k+1}, Y^{k+1}, Z^k)&\leq f(X^{k+1}, Y^k, Z^k)-\tau_2\|Y^{k+1}-Y^k\|_\textrm{F}^2,\\
f(X^{k+1}, Y^{k+1}, Z^{k+1})&\leq f(X^{k+1}, Y^{k+1}, Z^k)-\tau_3\|Z^{k+1}-Z^k\|_\textrm{F}^2.
\end{aligned}
\end{equation}
It follows that
\begin{equation}\label{fxyz-xyz}
\begin{aligned}
&f(X^{k+1}, Y^{k+1}, Z^{k+1})+\tau_1\|X^{k+1}-X^k\|^2_\textrm{F}+\tau_2\|Y^{k+1}-Y^k\|_\textrm{F}^2+\tau_3\|Z^{k+1}-Z^k\|_\textrm{F}^2\\
&\leq f(X^k, Y^k, Z^k).
\end{aligned}
\end{equation}
Therefore, it can be concluded that the update rule makes  the objective function nonincreasing strictly.\qed\\

\noindent{\textbf{Proof of Theorem \ref{Theor1} (b).}}
The boundedness of the sequence $\{(X^k, Y^k, Z^k)\}$ is proved by contradiction. On the one hand, suppose that the sequence $\{(X^k, Y^k, Z^k)\}$ is unbounded, and hence
\begin{equation}
\lim_{k\rightarrow\infty} \| (X^k, Y^k, Z^k)\|_\textrm{F}=\infty.
\end{equation}
It then follows from the coercive of $f(X,Y,Z)$ that the sequence $\{f(X^k, Y^k, Z^k)\}$ should diverge to infinity. Denote \begin{equation}
	\|E^{k+1}-E^k\|^2_\textrm{F}=\tau_1\|X^{k+1}-X^k\|^2_\textrm{F}+\tau_2\|Y^{k+1}-Y^k\|_\textrm{F}^2+\tau_3\|Z^{k+1}-Z^k\|_\textrm{F}^2.
	\end{equation} 
	On the other hand,
since \eqref{fxyz-xyz}, it has
\begin{equation}
\begin{aligned}
&f(X^{k+1}, Y^{k+1}, Z^{k+1})\\
&\leq f(X^{k+1}, Y^{k+1}, Z^{k+1})+\|E^{k+1}-E^k\|^2_\textrm{F}\leq f(X^k, Y^k, Z^k)\\
&\leq f(X^{k}, Y^{k}, Z^{k})+\|E^{k}-E^{k-1}\|^2_\textrm{F}\leq\cdots\leq f(X^0, Y^0, Z^0),
\end{aligned}
\end{equation}
which implies that $f(X^k, Y^k, Z^k)$ is finite for any $k$ and it leads to a contradiction. Therefore, the sequence $\{(X^k, Y^k, Z^k)\}$ is bounded.\qed\\

%

%
%
\noindent{\textbf{Proof of Theorem \ref{Theor1} (c).}} Let $K$ be a positive integer and $K>1$. Summing \eqref{fxyz-xyz} over $k=0,\ldots, K-1$ yields 

\begin{equation}
\begin{aligned}
&\sum_{k=0}^{K-1}\tau_1\|X^{k+1}-X^k\|_\textrm{F}^2+\tau_2\|Y^{k+1}-Y^k\|_\textrm{F}^2+\tau_3 \|Z^{k+1}-Z^k\|_\textrm{F}^2\\
&\leq\sum_{k=0}^{K-1}(f(X^k, Y^k, Z^k)-f(X^{k+1}, Y^{k+1}, Z^{k+1}))\\
&\leq f(X^0, Y^0, Z^0)- f(X^K, Y^K, Z^K)\\
&<+\infty,
\end{aligned}
\end{equation}
where the last inequality is due to  $f(X, Y, Z)$ being bounded from below. Hence,
\begin{equation}
\lim_{k\rightarrow\infty} \tau_1\|X^{k+1}-X^k\|_\textrm{F}^2+\tau_2\|Y^{k+1}-Y^k\|_\textrm{F}^2+\tau_3\|Z^{k+1}-Z^k\|_\textrm{F}^2=0,
\end{equation}
which suffices to 
\begin{equation}
\lim_{k\rightarrow\infty}\|(X^{k+1}, Y^{k+1}, Z^{k+1})-(X^k, Y^k, Z^k)\|_\textrm{F}=0.
\end{equation}

\noindent{\textbf{Proof of Theorem \ref{Theor1} (d).}}
\textbf{Step 1: A subgradient lower bound for the iterate gap.} Suppose that $\{(X^{k+1},Y^{k+1},Z^{k+1})\}$ is generated by Algorithm \ref{am1}. Following the first-order optimality conditions of \eqref{sub-x}, one has 
\begin{equation}
0\in\nabla_X f(X^{k+1},Y^{k+1},Z^{k+1})+\textrm{N}_{\mathcal{M}}(X^{k+1})+2\tau_1(X^{k+1}-X^k).
\end{equation}
Similarly, for  \eqref{sub-y} and \eqref{sub-z}, the following statements hold
\begin{equation}
\begin{aligned}
0&\in\nabla_Y f(X^{k+1},Y^{k+1},Z^{k+1})+\textrm{N}_{\mathcal{S}}(X^{k+1})+2\tau_2(Y^{k+1}-Y^k),\\
0&\in\nabla_Z f(X^{k+1},Y^{k+1},Z^{k+1})+\textrm{N}_{\mathcal{R}}(X^{k+1})+2\tau_3(Z^{k+1}-Z^k).
\end{aligned}
\end{equation}
There exist 
\begin{equation}
A^{k+1}=(A_X^{k+1},A_Y^{k+1},A_Z^{k+1})\in \nabla f(X^{k+1},Y^{k+1},Z^{k+1})+\textrm{N}_{\mathcal{M}\times \mathcal{S}\times \mathcal{R}}(X^{k+1},Y^{k+1},Z^{k+1})
\end{equation}
with 
\begin{equation}
	\begin{aligned}
A_X^{k+1}&\in\nabla_X f(X^{k+1},Y^{k+1},Z^{k+1})+\textrm{N}_{\mathcal{M}}(X^{k+1}),\\
A_Y^{k+1}&\in\nabla_Y f(X^{k+1},Y^{k+1},Z^{k+1})+\textrm{N}_{\mathcal{S}}(Y^{k+1}),\\
A_Z^{k+1}&\in\nabla_Z f(X^{k+1},Y^{k+1},Z^{k+1})+\textrm{N}_{\mathcal{R}}(Z^{k+1}).
	\end{aligned}
\end{equation}
such that
\begin{equation}
\begin{aligned}
0&=A_X^{k+1}+2\tau_1(X^{k+1}-X^k),\\
0&=A_Y^{k+1}+2\tau_2(Y^{k+1}-Y^k),\\
0&=A_Z^{k+1}+2\tau_3(Z^{k+1}-Z^k).
\end{aligned}
\end{equation}
Thus, for $\tau=2\max\{\tau_1,\tau_2,\tau_3\}$, one can get 
\begin{equation}
	\|(A_X^{k+1},A_Y^{k+1},A_Z^{k+1})\|_\textrm{F}\leq \tau\|(X^{k+1},Y^{k+1},Z^{k+1})-(X^k,Y^k,Z^k)\|_\textrm{F}.
\end{equation}

\textbf{Step 2: KL property.} It is well known that $\mathcal{M}$, $\mathcal{S}$ and $\mathcal{R}$ are semi-algebraic sets and their indicator functions are semi-algebraic. The quadratic functions $f(X,Y,Z)$ are also semi-algebraic. Using the fact that the composition of semi-algebraic functions is semi-algebraic, it derives that
\begin{equation}
f(X,Y,Z)+\delta_{\mathcal{M}}(X)+\delta_{\mathcal{S}}(Y)+\delta_{\mathcal{R}}(Z)
\end{equation}
is a semi-algebraic function, where $\delta_{\Omega}(\cdot)$ denotes the indicator function for set $\Omega$. Thus it satisfies the KL property at each point. 

Next, one can show that any accumulation point $(X^*,Y^*,Z^*)$ of the corresponding sequence $\{(X^k, Y^k, Z^k)\}$ is a stationary point of \eqref{dsco-2}. Let  $X^*$ be a stationary point of \eqref{penc-x}. According to \cite[Theorem 3.1]{xiao2022class}, for $\beta\geq \max\{2(\lambda_0+\lambda_1),2m\lambda_2\}$, it can easily verify that $X^*$ is also a stationary  point of \eqref{sub-x}. From \cite[Definition 2.1]{xiao2022class}, it derives 
\begin{equation}
\textrm{Tr}(U^\top\nabla_X f(X^*, Y^*, Z^*))\geq 0, ~{X^*}^\top X^*=I_m
\end{equation}
for any $U\in \textrm{T}_{\mathcal{M}}(X^*)$, which implies that 
\begin{equation}\label{KKTX}
0\in \nabla_X f(X^*, Y^*, Z^*)+\textrm{N}_{\mathcal{M}}(X^*).
\end{equation}
In addition, $Y^*\in \mathcal{S}$ is the closed-form solution of \eqref{sub-y}, which derives
\begin{equation}\label{KKTY}
0\in \nabla_Y f(X^*, Y^*, Z^*)+\textrm{N}_{\mathcal{S}}(Y^*).
\end{equation}
Similarly,  the following statement holds
\begin{equation}\label{KKTZ}
0\in \nabla_Zf(X^*, Y^*, Z^*)+\textrm{N}_{\mathcal{R}}(Z^*).
\end{equation}
Combining \eqref{KKTX}, \eqref{KKTY}, and \eqref{KKTZ}, it obtains 
\begin{equation}
0\in \nabla f(X^*, Y^*, Z^*)+\textrm{N}_{\mathcal{M}\times \mathcal{S}\times \mathcal{R}}(X^*,Y^*,Z^*).
\end{equation}
Thus, $(X^*, Y^*, Z^*)$ is a stationary point of \eqref{dsco-2}. According to \cite{bolte2014proximal}, combining \textbf{Step 1} and \textbf{Step 2},  it can be concluded that  $\{(X^k, Y^k, Z^k)\}$ generated by Algorithm \ref{am1} converges to a stationary point of \eqref{dsco-2}. \qed


%

\bibliographystyle{elsarticle-num}
\bibliography{mybibfile}

\end{document}